\input amstex
\documentstyle{amams} 
\document\annalsline{153}{2001}
\received{June 8, 1999}
\startingpage{329}
\input amssym.def
\input amssym.tex

\font\tenrm=cmr10

\def\joinrel{\mathrel{\mkern-4mu}}
\def\relbar{\mathrel{\smash-}}
\def\lrar{\relbar\joinrel\relbar\joinrel\rightarrow}
\def\llrar{\relbar\joinrel\relbar\joinrel\relbar\joinrel\relbar\joinrel\rightarrow}
\def\lllrar{\relbar\joinrel\relbar\joinrel\relbar\joinrel\relbar\joinrel\relbar\joinrel\relbar\joinrel\rightarrow}
\def\lorar{\relbar\joinrel\relbar\joinrel\relbar\joinrel\relbar\joinrel\relbar\joinrel\relbar\joinrel\relbar\joinrel\relbar\joinrel\relbar\joinrel\relbar\joinrel\rightarrow}
\def\lolrar{\relbar\joinrel\relbar\joinrel\relbar\joinrel\relbar\joinrel\relbar\joinrel\relbar\joinrel\relbar\joinrel\relbar\joinrel\relbar\joinrel\relbar\joinrel\relbar\joinrel\relbar\joinrel\relbar\joinrel\relbar\joinrel\relbar\joinrel\relbar\joinrel\rightarrow}
\catcode`\@=11
\font\twelvemsb=msbm10 scaled 1100

\font\ninemsb=msbm10 scaled 800
\newfam\msbfam
\textfont\msbfam=\twelvemsb  \scriptfont\msbfam=\ninemsb
  \scriptscriptfont\msbfam=\ninemsb
\def\msb@{\hexnumber@\msbfam}
\def\Bbb{\relax\ifmmode\let\next\Bbb@\else
 \def\next{\errmessage{Use \string\Bbb\space only in math
mode}}\fi\next}
\def\Bbb@#1{{\Bbb@@{#1}}}
\def\Bbb@@#1{\fam\msbfam#1}
\catcode`\@=12

 \catcode`\@=11
\font\twelveeuf=eufm10 scaled 1100
\font\teneuf=eufm10
\font\nineeuf=eufm7 scaled 1100
\newfam\euffam
\textfont\euffam=\twelveeuf  \scriptfont\euffam=\teneuf
  \scriptscriptfont\euffam=\nineeuf
\def\euf@{\hexnumber@\euffam}
\def\frak{\relax\ifmmode\let\next\frak@\else
 \def\next{\errmessage{Use \string\frak\space only in math
mode}}\fi\next}
\def\frak@#1{{\frak@@{#1}}}
\def\frak@@#1{\fam\euffam#1}
\catcode`\@=12


\def\ilim{\underset\leftarrow\to{\text{\rm lim}}}

\def\lefttensor{\buildrel \Bbb L\over \otimes}

\def\et{{\text{\rm  \'{e}t}}}

\def\Spec{\mathop{\roman{Spec}}\nolimits}
\def\Spf {\mathop{\roman{Spf}}\nolimits}
\def\Z{{\Bbb Z}}
\def\Fp{{\Bbb F_p}}

\def\Fpbar{\overline{\Bbb F}_p}

\def\Wx{W_x}

\def\Wbar{\overline{W}}
\def\iso{\buildrel \sim \over \longrightarrow}

\def\Gm{\Bbb{G}_m}
\def\Gmhat{\hat{\Bbb{G}}_m}

\def\E{\Cal E}

\def\L{\Cal L}
\def\F{\Bbb  F}

\def\M{\Cal M}

\def\O{\Cal O}
\def\C{\Cal C}

\def\det{\roman{det}}
\def\trace{\roman{trace}}
 
\font\emi= cmmi10 scaled 1700 



\title{Unit \hbox{\emi L}-functions and a \\ conjecture of Katz} 
\def\titleheadline#1{\def\one{#1}\ifx\one\empty\else
\gdef\thetitle{{\frenchspacing%
\let\\ \relax
{#1}}}\fi}
\newif\ifshort

\let\shorttitle\titleheadline
\shorttitle{ \eightsc\uppercase{Unit} {\eightpoint \it L}\/
\eightsc\uppercase{-functions} }
 \twoauthors{Matthew Emerton}{Mark Kisin}
 \institutions{University of Michigan, Ann Arbor, MI\\
{\eightpoint {\it E-mail address\/}: emerton\@math.lsa.umich.edu}\\
\vglue6pt
University of Sydney, Sydney, NSW, Australia\\
SFB 478, Westf\"alische Wilhelms Universit\"at, M\"unster, Germany\\
{\eightpoint {\it E-mail address\/}:  markk\@maths.usyd.edu.au}}



Let $p$ be a prime number, $\Fp$ the finite field of 
order $p$,
and $\Z_p$ the ring of $p$-adic integers.
Suppose that $X$ is a separated finite type $\Fp$-scheme
and that $\L$ is a lisse $\Z_p$-sheaf on the \'etale site of $X$. 
One defines 
in the usual way an $L$-function $L(X,\L)$ attached to $\L$. This
is a power series in a formal variable $T$, which by construction
is an element of $1+T\Z_p[[T]].$ If $f:X\rightarrow \Spec \Fp$ is the structural
morphism of $X$ then $f_!\L$ is a constructible
complex of\break $\Z_p$-sheaves on the \'etale
site of $\Spec \Fp$ of finite tor-dimension, and so we may form the $L$-function
$L(\Spec \Fp,f_!\L)$, which is also an element of $1+T\Z_p[[T]]$ (and in
fact a rational function).  The ratio $L(X,\L)/L(\Spec \Fp,f_!\L)$ thus
lies in $1+T\Z_p[[T]],$ and so may be regarded as a nowhere-zero analytic function
on the $p$-adic open unit disk $|T|<1.$ If $\L$ were a lisse $\Z_l$
sheaf with
$l\neq p$ then in fact this ratio would be identically 1 (this is
Grothendieck's approach to
the rationality of the zeta function). One does not have this in the setting
of $\Z_p$-sheaves. However, in this paper we prove the following result,
which was conjectured by Katz
([11, 6.1]):
\proclaim{Theorem} The ratio
$L(X,\L)/L(\Spec \Fp,f_!\L)$ extends to a nowhere\/{\rm -}\/zero function on the closed
unit disc $|T|\leq 1.$
\endproclaim

In particular this implies that the $L$-function for $L(X, \L)$ is 
$p$-adic meromorphic in the closed unit disc. According to Katz ([11, 6.1]), 
when $X = \Bbb A^n,$ this is an old result of Dwork's.

Katz also conjectured that $L(X,\L)$ extends to a meromorphic function
on the rigid-analytic affine $T$-line. The second part of Katz's conjecture 
was a generalization of a conjecture of Dwork predicting that $L(X, \L)$ was $p$-adic 
meromorphic when the unit $F$-crystal $\M$ associated to $\L$ was the unit part of 
an ordinary, overconvergent $F$-crystal. 
Wan [13] has found counterexamples
to the part of Katz's conjecture predicting that the $L$-function is meromorphic 
for general lisse sheaves. On the other hand, in a recent preprint  [14] 
he has 
proven that Dwork's more cautious prediction is actually true. Thus, while the 
question of when the $L$-function is meromorphic is by now quite well understood, 
much less is known about the location of the zeroes and poles. In fact the 
only prior results in this direction that we are aware of are due to
Crew [3], who proved Katz's conjecture in the special case when $X$ is
an affine curve and the sheaf $\L$ has abelian monodromy,
and Etesse and Le Stum [8], who proved Katz's conjecture under the 
strong assumption that $\L$ extends to some compactification of $X$.

In fact we prove a more general result showing that an analogue of Katz's conjecture 
is true even if one replaces $\Z_p$ by a complete noetherian local\break
$\Z_p$-algebra $\Lambda$ with finite residue field.
We refer to 
Corollary~1.8 for the precise formulation.

Our methods also have some 
applications to results on lifting representations 
of arithmetic fundamental groups. In particular, we show the following 

\vglue8pt

\proclaim{Theorem} Let $X$ be a smooth affine $\F_p$\/{\rm -}\/scheme{\rm ,} 
$\Lambda$ an artinian local $\Z_p$\/{\rm -}\/algebra{\rm ,} having finite 
residue field{\rm ,} and 
$$\rho: \pi_1(X) \rightarrow {\rm GL}_d(\Lambda)$$ 
a representation of the arithmetic {\rm \'{\it e}}tale fundamental group of $X.$ 
Consider a finite flat local $\Z_p$\/{\rm -}\/algebra $\tilde \Lambda$ 
and a surjection $\tilde \Lambda \rightarrow \Lambda.$ 
There exists a continuous lifting of $\rho$ 
$$\tilde \rho: \pi_1(X) \rightarrow {\rm GL}_d(\tilde \Lambda).$$ 

If $X$ is an open affine subset of $\Bbb P^1,$ 
then $\tilde \rho$ may be chosen so that the $L$\/{\rm -}\/function of the corresponding 
lisse sheaf of $\tilde \Lambda$\/{\rm -}\/modules is rational{\rm .}
\endproclaim

\vglue8pt

Let us now describe the contents of the paper in more detail.

In Section~1 we define the necessary 
$L$-functions, and state our main results. 

In Section~2 we give the proof of some standard results on behaviour of 
$L$-functions in exact triangles, and under stratification of the underlying space.
These are used later to reduce our calculations to a special case.

In Section~3 we introduce a key ingredient in our work, which 
is the relationship between locally constant or lisse \'etale sheaves and 
unit $F$-crystals. We show that if $\Lambda$ is an artinian (respectively
a complete noetherian) finite local $\Z_p$-algebra with finite residue field and 
if $X$ is a formally smooth $p$-adic formal scheme equipped with a lifting $F$ of the
absolute Frobenius of its special fibre,
then there is a one-to-one correspondence between locally constant \'etale sheaves of
finite free $\Lambda$-modules (respectively
lisse \'etale sheaves of $\Lambda$-modules) on $X$ and 
finite locally free $\Lambda\otimes_{\Z_p}\O_X$-modules $\E$ on $X$
equipped with a $\Lambda$-linear isomorphism 
$F^*\E \iso \E.$ This generalizes a result of Katz ([12, 4.1]),
which treats
the
cases when $\Lambda = \Z_p$ or $\Z/p^n$ for some positive integer $n$.

In Section~4 we begin our proof of Katz's conjecture.
For the usual reasons of technical convenience, 
we prove a more general result on complexes of $p$-adic constructible sheaves
of $\Lambda$-modules. By passing to the inverse limit, we see that it is enough to 
consider the case where $\Lambda$ is artinian.
Using standard techniques we reduce ourselves to having to verify Katz's conjecture
for a locally constant flat $\Lambda$-sheaf on an open subset of $\Gm^d$ 
with the additional property that the associated unit $(\Lambda,F)$-crystal $\E$ (as explained 
above) has an underlying
locally free sheaf of $\Lambda\otimes_{\Z_p}\O_X$-modules which is actually free.  
Under this assumption we may choose a surjection $\tilde \Lambda \rightarrow \Lambda,$ 
with $\tilde \Lambda$ a finite flat local $\Z_p$-algebra, and  lift this $\Lambda$-sheaf
to a lisse $\tilde \Lambda$-sheaf by lifting the associated unit $(\Lambda,F)$-crystal.

In this situation we prove Katz's conjecture directly,
by using an explicit trace calculation in the style of Dwork [5]. 
To reduce to this trace calculation we utilize a technique of Deligne
[4] in order to find an Artin-Schreier sequence for the extension by
zero of a lisse $\tilde \Lambda$-sheaf under an open immersion, which allows us to relate
$p$-adic \'etale cohomology with compact support to coherent cohomology of
formal $\Z_p$-schemes. A key part of this calculation is the subject of Section~5.

Finally, in Section~6 we show how our methods can be used to lift representations 
of arithmetic fundamental groups.

Finally, let us mention that the methods used in this paper are a special instance 
of a more general theory developed by the authors in the papers [6], [7]. 
This theory describes the \'etale cohomology of $p$-power torsion or\break $p$-adic
constructible sheaves on finite type $\Fp$-schemes in terms of quasi-coher\-ent
cohomology; more precisely, a certain Riemann-Hilbert
correspondence is established between the appropriate derived categories
of such \'etale sheaves, and a derived category of certain quasi-coherent
analogues of unit $F$-crystals.  However, we have written this
paper so as to be independent of these more general techniques.

\phantom{page break}

\demo{Acknowledgments}  The authors would like to thank Pierre Berthelot
for suggesting that they apply their techniques to the study of $L$-functions.
The first author would also like to thank Mike Roth for useful discussions.
 \enddemo

\section{Notation and statement of results}

(1.1)
Let $p$ be a prime number, fixed for the remainder of the paper,
$\Fp$ denote a finite field of order
$p$, and $\Fpbar$ be a choice of an algebraic closure of $\Fp$
(fixed once and for all).
For any positive integer $n$ we denote by $\Bbb F_{p^n}$ the unique degree $n$ 
extension of $\Fp$ contained in $\Fpbar$. 
We denote by $\sigma: \Fpbar \rightarrow \Fpbar$ the (arithmetic) Frobenius
automorphism $a \mapsto a^p$ of $\Fpbar,$ 
and denote by $\phi$ the geometric
Frobenius automorphism of $\Fpbar,$ that is, the inverse of $\sigma.$
Both $\sigma$ and $\phi$ are topological generators of the Galois group
of $\Fpbar$ over $\Fp$.

We denote by 
$\Wbar$ the ring of Witt vectors of $\Fpbar,$ and by  
$\sigma$ and $\phi$ respectively
the canonical Frobenius automorphism of $\Wbar$
and its inverse. If $n$ is any positive integer, we let $W(\Bbb F_{p^n})$
denote the Witt ring of $\Bbb F_{p^n}$, regarded as a subring of $\Wbar.$
As usual we will denote by $\Z_p$ the ring $W(\Fp)$ of $p$-adic numbers.
\medbreak
(1.2)
Let $X$ be a finite type scheme over $\Bbb F_p,$ and $\Lambda$ a 
noetherian ring, with $m\Lambda = 0$ for some positive integer $m.$ 
We denote by $D^-(X, \Lambda)$ the derived category of bounded-above 
complexes of
\'etale sheaves of $\Lambda$-modules on $X,$  and we denote by 
$D^b_{ctf}(X,\Lambda)$ the full subcategory of $D^-(X, \Lambda)$ consisting 
of complexes of $\Lambda$-sheaves which have 
finite tor-dimension, and which are constructible (in the sense that 
their cohomology sheaves are constructible\break $\Lambda$-modules [1]).
Given any object $\L^{\bullet}$ of $D^b_{ctf}(X,\Lambda)$, we may define its\break
$L$-function in the usual way, which we briefly recall.

If  $\L^{\bullet}$ is a constructible complex of $\Lambda$-sheaves of finite
tor-dimension on $X$ then by [4, Prop.-D\'ef.~4.6  (ii), p.~93], 
we may assume that $\L^{\bullet}$ is a finite length
complex of flat constructible $\Lambda$-sheaves. If $x$ is a closed point of
$X$ we let $\L^i_x$ denote the flat $\Lambda$-sheaf obtained by restricting
$\L^i$ to the point $x$. If we 
choose an identification of the residue field of $x$ with the subfield
$\Bbb F_{p^{d(x)}}$ of $\Fpbar,$ and denote by $\bar x$ the 
$\Fpbar$-valued point of $X$ lying over $x,$ 
then $\L^i_ {\bar x}$ is a free
$\Lambda$-module equipped with an action of the $d^{\rm th}$ power
$\phi^d$ of the geometric Frobenius automorphism. 
We let $|X|$ denotes the set of closed points of $X$, and define
$$ L(X,\L^{\bullet}) = \prod_{i}\prod_{x\in |X|}
\det_{\Lambda}(1-\phi^{d(x)}T^{d(x)},\L^i_{\bar x})^{(-1)^{i+1}}.$$
This is an element of $1+T\Lambda[[T]].$ 

From the definition we see that $L$-functions are compatible with 
change of ring: if $\Lambda, \Lambda'$ are noetherian rings killed by some 
positive integer, $g:\Lambda \rightarrow \Lambda'$ a map of rings, and $\L^\bullet$ 
is as above, then we have 
\smallbreak
\centerline{${\displaystyle L(X, \L^\bullet\overset \Bbb L \to \otimes \Lambda') = g(L(X, \L^\bullet)).}$}
\bigbreak

(1.3) Let $X$ be a finite type $\Bbb F_p$-scheme and let
$i: X_{\rm red} \rightarrow X$ be the closed immersion of the
underlying reduced subscheme of $X$ into $X$. Then the functors $i_!$ and $i^{-1}$
are quasi-inverse, and induce an equivalence of triangulated
categories between $D^b_{ctf}(X_{\rm red},\Lambda)$ and $D^b_{ctf}(X,\Lambda)$
(preserving the canonical $t$-structure of both categories).
If $\L^{\bullet}$ is any object of $D^b_{ctf}(X,\Z/p^n)$ then
$L(X,\L^{\bullet}) = L(X_{\rm red},i^{-1}\L^{\bullet}).$ 
\medbreak

(1.4)
If $f:X \rightarrow Y$ is a separated morphism of finite-type $\Bbb F_p$-schemes then push-forward with proper
supports induces a functor
$f_!: D^b_{ctf}(X,\Lambda) \rightarrow D^b_{ctf}(Y,\Lambda)$ ([4,
Th.\ 4.9, p.~95]).  

\specialnumber{1.5}
\proclaim{Theorem} 
Let $f:X \rightarrow Y$ be a separated morphism of finite\/{\rm -}\/type $\Bbb F_p$\/{\rm -}\/schemes{\rm ,} 
and $\Lambda$ a finite{\rm ,} local{\rm ,} artinian $\Z_p$-algebra{\rm ,} with maximal ideal $\goth m.$
If $\L^{\bullet}$ is any object of $D^b_{ctf}(X,\Lambda)$ then the ratio of
$L$\/{\rm -}\/functions 
 $$ L(X,\L^{\bullet})/L(Y,f_!\L^{\bullet}),$$
  a priori an element
of $1+T\Lambda[[T]] ,$ in fact lies in $1+\goth mT\Lambda [T].$
\endproclaim 

(1.6) Suppose that $\Lambda$ is finite and reduced of characteristic
$p$. 
Then from Theorem~1.5 we conclude that the $L$-functions of objects of $D^b_{ctf}(X_{\et},
\Lambda)$ are invariant under proper push-forward.
This was originally proved by Deligne ([4, Th.\ 2.2, p.~116]).

An interesting point is that Deligne gives
a counterexample involving a locally constant sheaf of 
free rank one $\F_p[X]/X^2$-modules, to show that his formula does not hold 
without the hypothesis that $\Lambda$ is reduced ([4, 4.5, p.~127]).
Thus if we take $\Lambda = \Bbb F_p[X]/X^2,$ 
then Deligne's counterexample shows that, in general, the quotient 
of the $L$-functions in Theorem~1.5 is not equal to $1.$ Nevertheless, 
our theorem asserts that it is equal to $1$ modulo the principal ideal $(X).$
This also 
follows from Deligne's theorem by ``specialization of $L$-functions'' (1.2),
and in fact this same argument shows that the ratio of $L$-functions
$L(X,\L^{\bullet})/L(Y,f_!\L^{\bullet})$ occurring in the statement of
Theorem~1.5 always lies in
$1+\goth mT\Lambda[[T]]$.
Thus the key result of Theorem~1.5 is that this ratio is in fact a polynomial.

\medbreak (1.7) We want to define $L$-functions for lisse sheaves, or more precisely, the 
lisse analogue of constructible sheaves. For this, suppose 
that we are given a noetherian ring $\Lambda$ and an ideal $I \subset \Lambda$,
such 
that $\Lambda$ is $I$-adically complete and $p$ is nilpotent in $\Lambda/I.$
We define the category $D^b_{I-sm}(X, \Lambda)$ to be the 
2-limit of the categories $D^b_{ctf}(X, \Lambda/I^n),$ 
$n = 1,2, \dots$. One has a formalism of $f_!$, $f^{-1}$ and $\lefttensor$
in $D^b_{I-sm}(X, \Lambda)$ 
(see [9] for details).

If $\L^{\bullet}$ is an object of $D^b_{I-sm}(X,\Lambda)$ then by
construction $\L^{\bullet}\lefttensor_{\Lambda}\Lambda/I^n$ belongs
to $D^b_{ctf}(X,\Lambda/I^n)$ for each positive integer $n$. If
$n \geq m$ then the compatibility of formation of $L$-functions with change
of rings shows that
$$L(X,\L^{\bullet}\lefttensor_{\Lambda}\Lambda/I^n) \equiv
L(X,\L^{\bullet}\lefttensor_{\Lambda}\Lambda/I^m)\pmod{I^m}.$$
Thus we may define $L(X,\L^{\bullet})\in 1+T\Lambda[[T]]$ to be the 
limit of the $L$-functions $L(X,\L^{\bullet}\lefttensor_{\Lambda}
\Lambda/I^n).$

Now suppose that $f:X\rightarrow Y$ is a separated map of finite type\break $\Fp$-schemes.
As before, we have a functor $f_!:D^b_{I-sm}(X,\Lambda) \rightarrow D^b_{I-sm}(Y,\Lambda).$
For each positive integer $n$ there is a canonical isomorphism 
$$(f_!\L^\bullet)\overset\Bbb L\to \otimes_{\Lambda}\Lambda/I^n \iso
f_!(\L^\bullet\lefttensor_{\Lambda}\Lambda/I^n),$$ so that 
$$L(Y,f_!\L^\bullet)\equiv L(Y,f_!(\L^\bullet\lefttensor_{\Lambda}\Lambda/I^n))\pmod{I^n}.$$

If we now suppose that $\Lambda$ is a complete local $\Z_p$-algebra with
finite residue field (so that in the above discussion $I=\goth m$ is
the maximal ideal of $\Lambda$), by letting
$n$ tend to infinity in the above discussion
we derive the following corollary of Theorem~1.5, which includes
Katz's conjecture as a special case (take $\Lambda$ to be $\Z_p$ and $\L^{\bullet}$ to
be a single lisse $p$-adic sheaf in degree zero):
 
\specialnumber{1.8}
\proclaim{Corollary} Let $\Lambda$ be a complete local $\Z_p$\/{\rm -}\/algebra{\rm ,} with 
finite residue field{\rm ,} and maximal ideal $\goth m.$
Let $f:X \rightarrow Y$ be a separated morphism of finite type $\Bbb F_p$\/{\rm -}\/schemes{\rm .}
If $\L^{\bullet}$ is any object of $D^b_{\goth m-sm}(X,\Lambda)$ then the ratio of\break
$L$\/{\rm -}\/functions $$L(X,\L^{\bullet})/L(Y,f_!\L^{\bullet}),$$
 a priori an element
of $1+T\Lambda[[T]],$ in fact lies in $1+\goth m T\Lambda\langle T \rangle,$ where 
$\Lambda\langle T\rangle$ denotes the $\goth m$\/{\rm -}\/adic
completion of the polynomial ring $\Lambda[T].$
\endproclaim 

\section{Preliminaries on $L$-functions}

(2.1)
Below, $\Lambda$ will be a noetherian $\Bbb Z/m\Bbb Z$-algebra for some 
integer $m.$

For any separated morphism $f:X \rightarrow Y$ of finite-type $\Bbb F_p$-schemes
and any object $\L^{\bullet}$ of $D^b_{ctf}(X,\Lambda)$, let us write
$$Q(f,\L^{\bullet}) = L(X,\L^{\bullet})/L(Y,f_!\L^{\bullet}).$$
We begin by recalling some standard tools for analyzing such ratios of\break $L$-functions.
Although the proofs are well-known and quite straightforward,
for the sake of completeness we recall them. 
 
\specialnumber{2.2}
\proclaim{Lemma} If $f:X \rightarrow Y$ is a separated morphism of finite\/{\rm -}\/type\break
$\Bbb F_p$\/{\rm -}\/schemes{\rm ,} then for any distinguished triangle
$$\cdots \rightarrow \L^{\bullet}_1 \rightarrow \L^{\bullet}_2 
\rightarrow \L^{\bullet}_3
\rightarrow \L^{\bullet}_1[1] \rightarrow \cdots $$ of objects of
$D^b_{ctf}(X,\Lambda),$
there is the equality
$$Q(f,\L^{\bullet}_2) = Q(f,\L^{\bullet}_1)Q(f,\L^{\bullet}_3).$$
\endproclaim

\demo{Proof} This follows 
immediately from the fact that $f_!$ takes distinguished
triangles to distinguished triangles, together with the multiplicativity of\break
$L$-functions of objects in a distinguished triangle.
\enddemo

If $f: X \rightarrow Y$ is a morphism of schemes and $y$ is a point of $Y$
we let $f_y: X_y \rightarrow y$ denote the fibre of $f$ over the point $y$.

\specialnumber{2.3}
\proclaim{Lemma} For any separated morphism $f:X \rightarrow Y$ of finite\/{\rm -}\/type\break
$\Bbb F_p$\/{\rm -}\/schemes
and any object $\L^{\bullet}$ of $D^b_{ctf}(X,\Lambda),$ there is the equality
$$Q(f,\L^{\bullet}) = \prod_{y\in |Y|}Q(f_y,\L^{\bullet}_{|X_y}).$$
\endproclaim

\demo{Proof}
By partitioning the elements of $|X|$ according to their image in $|Y|$, one
sees that
$$L(X,\L^{\bullet})=\prod_{y\in |Y|}L(X_y, \L^{\bullet}_{| X_y}).$$ By
the proper base-change theorem, for any point $y$ of $|Y|$, 
$$(f_!\L^{\bullet})_y = f_{y!}\L^{\bullet}_{|X_y},$$ and so
$$L(Y,f_!\L^{\bullet}) = \prod_{y\in |Y|} L(y,f_{y!}\L^{\bullet}_{|X_y}).$$
Thus $$Q(f,\L^{\bullet})  = \prod_{y\in |Y|}L(X_y, \L^{\bullet}_{| X_y})/
L(y,f_{y!}\L^{\bullet}_{|X_y}) 
= \prod_{y\in |Y|}Q(f_{y},\L^{\bullet}_{|X_y}),$$
proving the lemma.
\enddemo

\specialnumber{2.4}
\proclaim{Lemma} If $f:X \rightarrow Y$ and $g: Y \rightarrow Z$ are two
separated morphisms of finite\/{\rm -}\/type $\Bbb F_p$\/{\rm -}\/schemes and $\L^{\bullet}$ is any
object of $D^b_{ctf}(X,\Lambda)$ then there is the equality
$$Q(gf,\L^{\bullet}) = Q(g,f_!\L^{\bullet})Q(f,\L^{\bullet}).$$
\endproclaim

\demo{Proof}
We compute (beginning with the right-hand side)
$$Q(g,f_!\L^{\bullet})Q(f,\L^{\bullet})  = \frac{L(Y,f_!\L^{\bullet})}
{L(Z,g_!f_!\L^{\bullet})}\frac{L(X,\L^{\bullet})}{L(Y,f_!\L^{\bullet})}
 = \frac{L(X,\L^{\bullet})}{L(Z,(gf)_!\L^{\bullet})} = Q(gf,\L^{\bullet}).
$$
This proves the lemma.
\enddemo

\specialnumber{2.5}
\proclaim{Lemma} For any separated quasi\/{\rm -}\/finite morphism $f:X \rightarrow Y$
of finite\/{\rm -}\/type $\Bbb F_p$\/{\rm -}\/schemes and any object $\L^{\bullet}$ of
$D^b_{ctf}(X,\Lambda)${\rm ,} $Q(f,\L^{\bullet})=1${\rm .}
\endproclaim

\demo{Proof}
Lemma~2.3 reduces us to the situation in which $Y=\Spec \Bbb F_{p^d}$
is a point, so that
$X$ is either empty, in which case there is nothing to prove, or else
zero-dimensional. In this second case
we may replace $X$ by its underlying 
reduced subscheme (via (1.4)), and writing this as a disjoint union of points
we reduce to the case that $X$ is also a point, say $X = \Spec \Bbb F_{p^{d'}}$, with
$d$ dividing $d'$.
Using Lemma~2.2 one reduces to the
case that $\L^{\bullet}=\L$ is a single flat $\Lambda$-module equipped with
a $\Lambda$-linear $\phi^{d'}$ action. Then $f_!\L$ is simply the induced module
$$f_!\L = \Z/p^n[\phi^d]\otimes_{\Z/p^n[\phi^{d'}]}\L,$$ and the equality
of $L$-functions $L(X,\L)=L(Y,f_!\L)$ is an elementary calculation.
\enddemo

\specialnumber{2.6}
\proclaim{Lemma} Suppose that $f: X \rightarrow Y$ is separated morphism
of finite type $\Fp$\/{\rm -}\/schemes{\rm ,} and that $X = S_0 \coprod S_1 \coprod \cdots
\coprod S_n$ is a stratification of $X$ by locally closed subsets $S_i$ {\rm (}\/more
precisely{\rm ,}  each $S_i$ is closed in the union
$\coprod_{j=0}^i S_j${\rm ).} Each of the $S_i$ is given its
reduced\/{\rm /}\/induced scheme structure{\rm .} If $\L^{\bullet}$ is an object of 
$D^b_{ctf}(X,\Lambda)$ then let $\L^{\bullet}_{|S_i}$ 
denote the restriction of $\L^{\bullet}$ to
each of the locally closed subsets $S_i${\rm ,} and let $f_i:S_i \rightarrow Y$
denote the restriction of 
the morphism $f$ to each of the locally closed subsets $S_i${\rm .}  Then
$$Q(f,\L^{\bullet})= \prod_{i=0}^n Q(f_i,\L^{\bullet}_{|S_i}).$$
In particular{\rm ,} this applies if the stratification on $X$ is obtained by pulling 
back a stratification on $Y.$
\endproclaim

\demo{Proof}
This follows from Lemmas~2.2,~2.4 and~2.5, since immersions are quasi-finite.
\enddemo

\section{Etale sheaves and unit $F$-crystals}

(3.1) We want to explain a generalization of a result of Katz 
relating \'etale sheaves and unit $F$-crystals. Although what we 
do can be done somewhat more generally, we restrict ourselves 
to smooth $\Z/p^n$-schemes.

Assume we have a smooth scheme $X$ over 
$\Z/p^n,$ and let $\Lambda$ be a\break $\Z/p^n$-algebra. We 
assume that $X$ is equipped with an endomorphism $F$ 
lifting the absolute Frobenius on its reduced subscheme. 

A unit $(\Lambda, F)$-crystal 
on $X$ is a sheaf $\E$ of finite 
locally free $\Lambda\otimes_{\Z/p^n}\O_X$ 
modules, equipped with a $\Lambda\otimes_{\Z/{p^n}}\O_X$-linear isomorphism 
$F^*\E \iso \E.$ Note that if the map $\Z/{p^n} \rightarrow \Lambda$ 
factors through $\Z/p^{n'}$ for some positive integer $n'< n,$ then 
the notion of a unit $(\Lambda,F)$-crystal depends only on the 
reduction of $X$ modulo $p^{n'}.$

\specialnumber{3.2}
\proclaim {Proposition} Suppose that $\Lambda$ is noetherian{\rm ,} local{\rm ,} and 
finite over $\Z/p^n.$ 
There is an equivalence of categories 
{\rm (}\/explicitly described below\/{\rm )} between the category of 
locally constant {\rm \'{\it e}}tale sheaves of finite free $\Lambda$\/{\rm -}\/modules 
on $X$ and the category of $(\Lambda, F)$\/{\rm -}\/crystals on $X.$
\endproclaim

\demo{Proof} Let $\L$ be a locally constant 
\'etale sheaf of finite free $\Lambda$-modules on $X$ 
(which is, of course equivalent to the data of such a sheaf 
on the reduction of $X$ modulo $p$).
We associate to $\L$ a $(\Lambda,F)$-crystal on $X$
as follows. 
Consider the \'etale sheaf $\E_{\et} =: \L\otimes_{\Z/p^n}\O_{X_{\et}}.$ 
This is a coherent sheaf of $\Lambda\otimes_{\Z/p^n}\O_{X_{\et}}$-modules,
which is in fact locally free of finite rank (since $\L$ is locally
free of finite rank over $\Lambda$). 
Regarding $\E_{\et}$ as a coherent sheaf on $\O_{X_{\et}},$ we see by 
[12, 4.1]  that $\E_{\et}$ comes by pull-back from a Zariski coherent 
$\O_X$-module $\E.$ Moreover, $\E$ is equipped with an isomorphism 
$F^* \E \iso \E.$ If we denote by $\Phi$ the composite 
$\E \rightarrow F^*\E \iso \E,$ then the induced map $\Phi_{\et}$ on $\E_{\et}$ 
has fixed subsheaf equal to $\L.$ 
As the formation of $\E$ is functorial in $\E_{\et},$ 
we see that $\E$ is equipped with the structure of a $\Lambda$-module, 
and that the isomorphism $F^*\E \iso \E$ is $\Lambda$-linear. 

To see that 
this gives $\E$ the structure of a unit $(\Lambda, F)$-crystal we have 
to check that $\E$ is locally free as a $\Lambda\otimes_{\Z/p^n}\O_X$-module. 
If $x \in X$ denote by $\tilde x$ an \'etale point lying over $x.$
Now
$\E_{\et,\tilde x}=\E_x\otimes_{(\Lambda\otimes_{\Z/p^n}\O_{X,x})}
(\Lambda\otimes_{\Z/p^n}\O_{X_{\et},\tilde x})$
is certainly a free $\Lambda\otimes_{\Z/{p^n}}\O_{X_{\et},\tilde x}$ 
module. Since $\Lambda\otimes_{\Z/{p^n}}\O_{X_{\et},\tilde x}$ is faithfully
flat over $\Lambda\otimes_{\Z/p^n}\O_{X,x}$ ($\O_{X_{\et},\tilde x}$ being
faithfully flat over $\O_{X,x}$), we conclude by descent that
$\E_x$ is locally free of constant rank over $\Lambda\otimes_{\Z/p^n}\O_{X,x}$.
Since $\Lambda$ is finite over $\Z/p^n$ we see that
$\Lambda\otimes_{\Z/p^n}\O_X$ is finite over the local ring $\O_{X,x}$,
and so is semi-local. Thus the freeness of ${\Cal E}_x$ over
$\Lambda\otimes_{\Z/p^n}\O_X$ follows from the   (simple) commutative
algebra fact: a finite module $M$ over a semi-local ring $A$ which is locally
free of constant rank is in fact free. (To see this, observe that $M/{\rm rad}(A)$
is locally free of constant rank over the direct sum of fields $A/{\rm rad}(A)$,
and so is certainly free. Lifting generators,
Nakayama implies that $M$ is itself
free.)

Next we construct the quasi-inverse functor.
Given a unit $(\Lambda,F)$-crystal $\E$ of rank $m$,
we pull it and its  Frobenius endomorphism  
$\Phi$ back to the \'etale site to get a coherent locally free 
$\Lambda\otimes_{\Z/p^n}\O_{X_{\et}}$-module $\E_{\et},$ equipped with 
a $\Lambda$-linear endomorphism $\Phi_{\et}.$ 
We have to 
show that if $\L = \ker (1-\Phi_{\et}),$ then the natural map 
$\L\otimes_{\Z/p^n}\O_{X_{\et}} \rightarrow \E_{\et}$ is an isomorphism. Indeed,
once 
we have this, then $\L$ is necessarily locally free over $\Lambda,$ by flat descent, 
as $\O_{X_{\et}}$ is flat over $\Z/p^n,$ and $\E_{\et}$ is locally free over
$\Lambda\otimes_{\Z/p^n}\O_{X_{\et}}.$

For this, suppose first that $\Lambda$ is flat (and necessarily finite) over $\Z/p^n.$
In this case we may regard $\E$ as a unit $(\Bbb Z/p^n, F)$-crystal, so  that the map 
above is an isomorphism by Katz's theorem [12, 4.1].

In general, we can write $\Lambda$ as a quotient of 
a finite flat local $\Z_p$-algebra
$\tilde \Lambda$; as $\Lambda$ is finite over
$\Z_p,$ there is a surjection $h:\Z_p[x_1, \dots, x_r] \rightarrow \Lambda,$ 
and for $i = 1, \dots ,r$ there exists a monic polynomial $p_i$ with 
coefficients in $\Z_p$ such that $h(p_i(x_i)) = 0.$ Thus 
$\Lambda$ is a quotient of  $\Z_p[x_1, \dots, x_r]/(p_1(x_1), \dots, p_r(x_r)),$
which is a finite flat $\Z_p$-algebra, hence, in particular, semi-local. 
Thus we may take $\tilde \Lambda$ to be a localization of 
$\Z_p[x_1, \dots, x_r]/(p_1(x_1), \dots, p_r(x_r))$ at a suitable maximal ideal. 
Then $\tilde \Lambda$ is finite flat over $\Z_p,$ as $\Z_p$ is a complete local ring 
(so localization does not destroy finiteness).

To show that   $\L\otimes_{\Z/p^n}\O_{X_{\et}} \iso \E_{\et}$
we may work locally. As $\E$ is a 
locally free
$\Lambda\otimes_{\Z/p^n}\O_X$-module, after localizing on $X,$ we may assume 
that $\E$ is free over $\Lambda\otimes_{\Z/p^n}\O_X,$ say
$\E = (\Lambda\otimes_{\Z/p^n}\O_X)^m.$ Write
$\tilde\E = (\tilde\Lambda/p^n \otimes_{\Z/p^n}\O_X)^m.$ 
In this case we can lift the isomorphism $F^*\E \iso \E$ of free 
$\Lambda\otimes_{\Z/p^n}\O_X$-modules to a morphism  
$ F^*\tilde\E \rightarrow \tilde\E$ of free
$\tilde\Lambda/p^n\otimes_{\Z/p^n}\O_X$-modules,
and any such lift 
is an isomorphism, by Nakayama's lemma.
This gives $\tilde \E$ the structure of a 
unit $(\tilde \Lambda/p^n, F)$-crystal.
We denote by $\tilde\Phi:\tilde \E \rightarrow \tilde \E$
the induced $\tilde\Lambda$-linear,
$F^*$-semi-linear endomorphism of $\tilde{\E}$, and
by $\E_{\et}, \tilde\E_{\et}, \Phi_{\et}, \tilde\Phi_{\et}$ the pull-backs to the \'etale 
site of $X$ of $\E, \tilde\E, \Phi, \tilde\Phi.$ If $x$
is any \'etale point of $X_{\et},$ 
we obtain a commutative diagram 
$$\matrix 
0&\lrar & \ker (1- \tilde\Phi_{\et})_{x} &\lrar&
\tilde\E_{x}&
 \buildrel{1-\tilde\Phi_{\et}}\over{\llrar}& \tilde \E_{x} &\lrar & 0 \\
&&\Big\vert&&\Big\vert&&\Big\vert\\
\noalign{\vskip-8pt}
&&\downarrow &&\downarrow&&\downarrow\\
0 &\lrar & \ker (1- \Phi_{\et})_{x} &\lrar&
\E_{x}&
 \buildrel{1-\Phi_{\et}}\over{\llrar} &   \E_{x} &\lrar & 0 .\endmatrix$$ 
Here the bottom row is obtained by applying
$\otimes_{\tilde\Lambda/p^n}\Lambda$ 
to the top row. By the previous discussion, applied to the finite flat $\Z/p^n$ 
algebra $\tilde \Lambda/p^n,$ we know that the top row is exact, hence the bottom one is also, as 
$\tilde\E_x$ is flat over $\tilde\Lambda/p^n,$ $\O_X$ being flat over $\Z/p^n.$ 
Thus, letting $\tilde \L$ denote $\ker(1-\tilde\Phi_{\et}),$ we see that the map
$\L\otimes_{\Z/p^n}\O_{X_{\et}} \rightarrow \E_{\et}$ is obtained 
by applying $\otimes_{\tilde \Lambda/p^n}\Lambda$ to the isomorphism 
$\tilde \L\otimes_{\Z/p^n}\O_{X_{\et}} \iso \E_{\et},$ and so is an isomorphism.
\enddemo

(3.3) We may formulate a version of the above theorem 
for \'etale sheaves on a smooth $\Bbb F_p$-scheme $X_0.$ 
We may define 
a unit $(\Lambda, F)$-crystal $\E$ to be a locally free 
sheaf of finite rank $\Lambda\otimes_{\Z_p} \O_{X_0,{\rm crys}}$-modules on the 
crystalline site of $X_0/\Z/p^n$ 
equipped with an isomorphism $F^*\E \iso \E$ (here $F$ denotes the 
Frobenius on the crystalline site, and $\O_{X_0,{\rm crys}}$ is the structure
sheaf on the crystalline site). If there exists a smooth $\Z/p^n$ scheme $X$ 
whose special fibre is $X_0$ and a lift $F$ of the absolute Frobenius on $X,$ 
then the two 
notions are equivalent (exercise (!), but see also [7]),
so that  our new definition is 
consistent with the previous one. Thus even if no global lift 
exists we get an equivalence of categories between locally constant \'etale 
sheaves of free $\Lambda$-modules and  unit $(\Lambda, F)$-crystals. 
Indeed, to construct this equivalence, we may work locally, and then we may 
assume ([10, III]) that $X_0$ lifts to a smooth $\Z/p^n$ scheme, 
equipped with a lift of Frobenius. Now we can appeal to our previous 
results.

The referee has also remarked that there is a connection between the
preceding theorem and the results of Berthelot and Messing which relate
finite \'etale group schemes and Dieudonn\'e crystals (see \S 2 of [2]).

\medbreak (3.4) It will be convenient to have an analogue of Proposition~3.2 for lisse 
sheaves. For this let $X$ be a formally smooth $p$-adic formal scheme over $\Z_p,$ 
equipped with a lift $F$ of the absolute Frobenius on its reduced subscheme, 
and $\Lambda$ a finite flat $\Z_p$-algebra. A unit $(\Lambda,F)$-crystal $\E$ 
is defined as in (3.1). Namely, it is a locally free, coherent $\Lambda\otimes_{\Z_p}\O_X$-module 
$\E$ equipped with an isomorphism $F^* \E \iso \E.$ Then we have:\enddemo

\specialnumber{3.5}
\proclaim {{C}orollary} Suppose that $\Lambda$ is local{\rm ,} and finite flat over $\Z_p.$ 
There is an equivalence of categories between $p${\rm -}\/adic lisse 
sheaves of $\Lambda$\/{\rm -}\/modules on $X$ and unit $(\Lambda, F)$\/{\rm -}\/crystals on
$X${\rm .}
\endproclaim

\demo {Proof} This follows immediately from Proposition~(3.3) 
once we note that 
a unit $(\Lambda,F)$-crystal $\E$ satisfies 
$\E \iso \ilim\, \E/p^n\E,$ because $\Lambda$ is finite over $\Z_p.$

Explicitly, if $\L$ is a $p$-adic lisse sheaf of $\Lambda$-modules, 
write $\L_n = \L/p^n\L.$ Then $\L_n\otimes_{\Z_p}\O_{X_{\et}}$ descends 
to a coherent, locally free Zariski sheaf of $\Lambda/p^n\otimes_{\Z_p}\O_X$-modules $\E_n$ on $\O_X,$ 
equipped with an isomorphism $F^*\E_n \iso \E_n.$ Then we may attach 
$\E = \ilim \E_n$ to $\L.$ It is locally free over $\Lambda\otimes_{\Z_p}\O_X$ 
as each of the $\E_n$ is locally free over $\Lambda/p^n\otimes_{\Z_p}\O_X.$ 
(In fact if $U \subset X$ is an open formal subscheme
then $\E$ becomes free over $U$ as 
soon as $\E_1$ is free over $U$.)

Conversely, if $\E$ is a unit $(\Lambda,F)$-crystal, then we may attach a locally constant 
\'etale sheaf of free $\Lambda/p^n$-modules $\L_n$ to $\E/p^n\E,$ and set $\L = \ilim \L_n.$

If $\E_{\et}$ is the pull-back of $\E$ to $X_{\et},$ we sometimes abuse notation, and 
write $\E_{\et} = \L\otimes_{\Z_p} \O_{X_{\et}}.$
\enddemo

\section{Proof of Theorem~1.5}

\vglue6pt
(4.1) In this section we present the proof of Theorem~1.5.
We denote by $\Lambda$ a finite, local, artinian $\Z_p$-algebra.
We begin with some preliminaries on formal schemes, and liftings of Frobenius.

\bigbreak (4.2)
Let $d$ be a positive integer. Let $\Bbb P^d$ denote $d$-dimensional 
projective space over $\F_p,$ equipped with homogeneous
coordinates $Z_0,\ldots, Z_d.$
Write\break $z_i:=Z_i/Z_0$ for
the affine coordinates corresponding to our choice of homogeneous coordinates;
then 
$$\Gm^d=\Spec \Fp[z_1,z_1^{-1},\ldots, z_d,z_d^{-1}] \subset \Bbb P^d.$$

We denote by $\hat{\Bbb P}^d$ the
formal, $d$-dimensional projective space over $\Z_p$ 
  obtained by completing the $d$-dimensional projective space $\Bbb P_{\Z_p}^d$ over $\Z_p$
along the special fibre $p=0$.
The underlying topological spaces of $\hat{\Bbb P}^d$ and $\Bbb P^d$ are equal.
We denote by $\O_{\hat{\Bbb P}^d}$ the structure sheaf of $\hat{\Bbb P}^d$,
thought of as a sheaf on the underlying topological space of $\Bbb P^d$,
and if $U$ is any open subset of $\Bbb P^d$, we write 
$\O_{\hat U} = \O_{\hat{\Bbb P}^d}|_U.$ 
We will sometimes write $\hat{U}$ to
denote $U$ considered with the formal scheme structure given by $\O_{\hat U}.$ 
We will employ this notation  especially in
the case that $U$ is an open subset of $\Gm^d \subset \Bbb P^d.$ 
\bigbreak

(4.3)
Let $X$ be of finite type over $\Fp$.  If $x$ is a closed point
of such a scheme we let $\kappa(x)$ denote the residue field of $x$. This is
a finite field of some degree $d(x)$ over $\Fp.$ We denote the ring of Witt
vectors of $\kappa(x)$ by $\Wx.$ We again use $\sigma$ (respectively $\phi$)
to denote the Frobenius automorphism of $\Wx$ (respectively its inverse automorphism).
If we choose an isomorphism of $\kappa(x)$
with $\Bbb F_{p^{d(x)}}$, this induces an isomorphism of $\Wx$ with
$W(\Bbb F_{p^{d(x)}})$, which is equivariant with respect to the automorphism
$\sigma$ of each of these rings.

Consider the endomorphism $F$ of $\hat{\Bbb P}^d$ given by $Z_i\mapsto Z_i^p.$ 
This is a lift of the absolute Frobenius to $\hat{\Bbb P}^d,$
which induces an endomorphism of each formal open $\hat U \subset \hat{\Bbb P}^d,$ 
and in particular of the formal $d$-dimensional multiplicative group 
$$\Gmhat^d = \Spf(\Z_p\langle z_1,z_1^{-1},\ldots, z_d,z_d^{-1}\rangle) \subset \hat{\Bbb P}^d,$$ 
where it is given by $F^*(z_i) = z_i^p.$ 
(Here $\Z_p\langle z_1,z_1^{-1},\ldots, z_d,z_d^{-1}\rangle$ denotes the
$p$-adic completion of the $\Z_p$-algebra $\Z_p[z_1,z_1^{-1},\ldots,
z_d,z_d^{-1}]$.)

If $k$ is a finite 
field, and $x$ is a $k$-valued point of $\Gmhat^d$
corresponding to the morphism
$$\Fp[z_1,z_1^{-1},\ldots,z_d,z_d^{-1}] \rightarrow k$$ of $\Fp$-algebras,
then $x$ lifts to a morphism
$$\Z_p\langle z_1,z_1^{-1},\ldots,z_d,z_d^{-1}\rangle \rightarrow W(k)$$
(the Teichm\"uller lifting), characterized by the property that the
diagram
$$\matrix \noalign{\vskip4pt}
\Z_p\langle z_1,z_1^{-1},\ldots,z_d,z_d^{-1}\rangle &\lrar&  
W(k)\\ \noalign{\vskip4pt}
\Big\vert{\scriptstyle F^*}&&\Big\vert{\scriptstyle\sigma}\\
\noalign{\vskip-8pt}
\hskip3pt\downarrow\phantom{\scriptstyle F^*}&&\hskip3pt\downarrow\phantom{{\scriptstyle\sigma}}\\ \noalign{\vskip4pt}
\Z_p\langle z_1,z_1^{-1},\ldots,z_d,z_d^{-1}\rangle &\lrar &W(k)\\ \noalign{\vskip4pt}\endmatrix$$
commutes. We denote this $W(k)$-valued point of $\Gmhat^d$ by $\tilde{x}.$

In particular this applies if $x$ is a closed point of $\Bbb G^d_m$ and we take 
$k$ to be $\kappa(x),$ the residue field at $x,$ so that $x$ is naturally
a $\kappa(x)$-valued point of $\Gm^d$.
If we choose an embedding of $\kappa(x)$ into a subfield $\Bbb F_{p^n}$
of $\Fpbar$ (for some positive integer $n$) then we may regard $x$ as
an $\Bbb F_{p^n}$-valued point of $\Gm^d,$ and the corresponding embedding of
$W_x$ as a subring of $W(\Bbb F_{p^n})$ realizes $\tilde{x}$ as a
$W(\Bbb F_{p^n})$-valued point of $\Gmhat^d.$  
\bigbreak

(4.4) We now begin the proof of Theorem~1.5.
Consider an arbitrary separated morphism $f:X \rightarrow Y$
of finite type $\Fp$-schemes.  We may find a stratification of $Y$ by
locally closed affine schemes $T_i$. If $S_i=f^{-1}(T_i)$
then Lemma~2.6 shows that it suffices
to prove Theorem~1.5 for each of the morphisms $f_i: S_i \rightarrow T_i$
obtained by restricting $f$ to the subschemes $S_i$. 

Since each $T_i$ is affine, the structural morphisms $g_i: T_i\rightarrow \Spec
\Fp$ and $h_i=f_ig_i:S_i\rightarrow \Fp$ are separated, and Lemma~2.4 shows that
for any object $\L^{\bullet}$ of $D^b_{ctf}(S_i,\Lambda)$,
$$Q(f_i,\L^{\bullet}) = Q(h_i,\L^{\bullet})/Q(g_i,f_{i!}\L^{\bullet}).$$ Thus it
suffices to prove Theorem~1.5 for the morphisms $g_i$ and $h_i$. 

\bigbreak (4.5) The preceding section shows that 
we are reduced to considering Theorem~1.5 in the case that $f:X\rightarrow \Spec
\Fp$ is the structural morphism of a finite type separated $\Fp$-scheme.
We will prove Theorem~1.5 by induction
on the dimension of $X$. Thus we assume that the result holds for all separated
$\Fp$-schemes of finite type and of dimension less than that of $X$.

Nothing is changed if we replace $X$ by its reduced subscheme. 
We find a dense open subscheme $U$ of $X$ such that each connected component
of $U$ admits a quasi-finite and dominant
morphism to $\Gm^d$ for some natural number $d$ (which
is necessarily the dimension of this connected component,
and so is less than or equal to the
dimension of $X$). To find such a $U$ we choose
a dense smooth open subscheme of $X$. Each connected component
of this smooth open subscheme is irreducible,
and contains a nonempty (and so dense) open subset
which admits an \'etale map to some $\Gm^d$. We take $U$ be the union of 
these open subsets. Since $X\setminus U$ has lower dimension than $X$,
Lemma~2.6 together with our inductive hypothesis
shows that it suffices to prove Theorem~1.5 for each connected component 
of $U$. Fix one such connected component $V$,  say, equipped
with a quasi-finite dominant morphism $g:V\rightarrow
\Gm^d$. Let\break $h:\Gm^d \rightarrow \Spec \Fp$ be the structural morphism of $\Gm^d$,
and $k=gh: V \rightarrow \Spec \Fp$ the structural morphism of $V$.
For any object $\L^{\bullet}$ of $D^b_{ctf}(V,\Lambda)$, Lemmas~2.4 and~2.5
show that 
$$Q(k,\L^{\bullet}) = Q(h,g_!\L^{\bullet})Q(g,\L^{\bullet}) = Q(h,g_!\L^{\bullet}).$$
So it suffices to prove Theorem~1.5 when $X = \Bbb G^d_m,$ and $Y = \Spec( \Bbb F_p).$
\bigbreak

(4.6)
Let $\L^{\bullet}$ be any object of $D^b_{ctf}(\Gm^d,\Lambda)$.
As observed in (1.2), we may assume that each of the finitely many nonzero
$\L^i$ is a constructible \'etale sheaf of flat $\Lambda$-modules.  Lemma~2.2 allows us to verify
Theorem~1.5 for each of the $\L^i$ separately. Thus we assume that $\L^{\bullet}
=\L$ is a single constructible \'etale sheaf of flat $\Lambda$-modules.

We may find a nonempty (and hence dense) open affine subscheme $U$ of $\Gm^d$ such that
$\L$ is locally constant when restricted to $U.$ 
We form the tensor product
$\E=\L_{|U}\otimes_{\Z_p}\O_{\hat U}.$ Note that $\E$ is a coherent 
$\O_{\hat U}/p^n$-sheaf for sufficiently large integers $n.$  Hence it descends 
to a coherent Zariski sheaf of locally free $\Lambda \otimes_{\Z_p} \O_{\hat U}$-modules.
In particular, $U$ contains a nonempty (and hence dense) open subset $V$ over 
which $\E$ becomes free as a $\Lambda\otimes_{\Z_p}\O_{\hat U}$-module. 
We may find a global section $a$ of $\O_{\Gm^d}$
such that the zero-set $Z(a)$ of $a$ is a proper subset of
$\Gm^d$
containing the complement of $V$.
Writing $D(a)=\Gm^d\setminus Z(a)\subset V,$ we see that $D(a)$ is nonempty
and that $\E$ is free over $\Lambda\otimes_{\Z_p}\O_{\hat U},$ when restricted to $D(a).$ 
Since $Z(a)$ has dimension less than $d,$ our induction hypothesis, together with Lemma~2.6,  
shows that it suffices to prove that $Q(f,\L_{|D(a)})$ lies in $1+\goth m T\Lambda[T].$ 
\bigbreak

(4.7) Since we will only be dealing with the locally constant sheaf
$\L_{|D(a)}$ from now on, we write simply $\L$ and $\E$ rather than $\L_{|D(a)}$ and
$\E_{|D(a)}$, so that $\E=\L\otimes_{\Z_p}\O_{\hat{D}(a)}.$

Recall 
the endomorphism $F$ of $\hat{D}(a)$ defined in Section~(4.3) by restricting
the endomorphism $F$ of $\hat{\Bbb P}^d.$
The tensor product of the identity morphism
on $\L$ and the endomorphism $1\otimes F^*$ of $\Lambda\otimes_{\Z_p}\O_{\hat{D}(a)}$ 
induces an $1\otimes F^*$-linear endomorphism of $\E$,
which we denote by $\Phi$, and endows $\E$ with the structure of a unit $(\Lambda, F)$-crystal,
from which $\L$ may be recovered as the \'etale subsheaf of $\Phi$-invariants 
(see \S 3). For convenience we will often abbreviate $1\otimes F^*$ to~$F^*.$

The global section $a$ is an element of
$\Fp [z_1,z_1^{-1},\ldots,z_d,z_d^{-1}]$.
Let $\tilde{a}$ denote a lift of $a$ to a global section of $\O_{\Gmhat^d}.$
Then {\it a priori} $\tilde{a}$ is an element of $\Z_p\langle
z_1,z_1^{-1},\ldots,z_d,z_d^{-1}\rangle,$ but we may and do choose $\tilde{a}$
to be an element of
$\Z_p[z_1,z_1^{-1},\ldots,z_d,z_d^{-1}]$ (since our only concern
is that it reduce modulo $p$ to~$a$).  

Choosing
a basis of $\E$ induces an isomorphism $\Lambda\otimes_{\Z_p}\O_{\hat D(a)}^m \iso \E.$ 
With respect to this
basis the $F^*$-linear endomorphism $\Phi$ may be written in the form 
$\Phi=(r_{ij})\circ F^*,$ where $(r_{ij})$ is an invertible
$m\times m$ matrix of elements
of $\Lambda \otimes_{\Z_p}\Z_p[z_1,z_1^{-1},\ldots, z_d,z_d^{-1},\tilde a^{-1}].$ 

As in the proof of Proposition~3.2, there exists a finite flat local\break $\Z_p$-algebra 
$\tilde \Lambda,$ and a surjection $\tilde\Lambda \rightarrow \Lambda.$ We write 
$\tilde \goth m$ for the maximal ideal of $\tilde \Lambda.$
Now we may lift each of the rational functions $r_{ij}$ to an element $\tilde r_{ij}$
of $\tilde \Lambda [z_1,z_1^{-1},\ldots,z_d,z_d^{-1},\tilde{a}^{-1}].$ 
Using these lifts
we may form the following unit $(\tilde\Lambda, F)$-crystal on the open subset $\hat{D}(a)$ of
$\Gmhat^d$:
$$ \tilde{\Phi}:
\tilde\Lambda\otimes_{\Z_p}\O_{\hat{D}(a)}^m \buildrel{(\tilde{r}_{ij})\circ F^*} \over\lllrar 
\tilde\Lambda\otimes_{\Z_p}\O_{\hat{D}(a)}^m,$$ 
which reduces to $\E$ after $\otimes_{\tilde\Lambda}\Lambda.$ Denote this 
unit $(\tilde\Lambda, F)$-crystal by $\tilde{\E}$, and let $\tilde{\L}$
denote the corresponding lisse \'etale $\Z_p$-sheaf on $\hat{D}(a)$
obtained as the\break $\tilde{\Phi}$-fixed 
\'etale subsheaf of the \'etale sheaf induced by $\tilde {\E}.$ We saw in
Section~3 that
$\tilde\L$ is a lisse sheaf of $\Lambda$-modules of rank $m.$ 
We will prove that the ratio
$L(D(a),\tilde{\L})/L(\Spec \Fp,g_!\tilde{\L})$, 
belongs to $1+\tilde \goth m T\tilde\Lambda\langle T\rangle.$ Specializing 
via the map $\tilde \Lambda \rightarrow \Lambda,$ 
we conclude
that $$Q(f,\L)\in 1 + \goth m T\Lambda/p^n[T],$$ completing the proof
of Theorem~1.5.\bigbreak

(4.8) In this section we make a series of changes of basis of the
unit\break $F$-crystal ${\Cal E}$ and its lift $\tilde{\E}$ in order to ensure that
the matrices $(r_{ij})$ and $(\tilde{r}_{ij})$ have certain properties.
Slightly abusing notation, we will again denote by $\tilde a$ the image of 
$\tilde a$ in $\Lambda[z_1, z_1^{-1}, \dots, z_d, z_d^{-1}].$
We begin by observing that
$$\align F^*\tilde{a} & \equiv \tilde{a}^p \pmod p\\
& \equiv 0 \pmod {\tilde{a}^p,p}.\endalign$$
Thus, if $n$ is an integer such that $p^n = 0$ in $\Lambda,$ then for some sufficiently 
large natural number $N$ we have 
$$(F^*\tilde{a})^N \equiv 0 \pmod{\tilde{a}^{N+1},p^n},$$ and so 
$$(F^*\tilde a/\tilde a)^N \in \tilde a\Lambda [z_1,z_1^{-1},\ldots,z_d,z_d^{-1}].$$
If we let $M$ be the maximal power
of $\tilde a$ occurring in the denominators of the elements $r_{ij}$
of $\Lambda[z_1,z_1^{-1},\ldots,z_d,z_d^{-1},\tilde a^{-1}]$ then
we see that
$$(F^*\tilde a/\tilde a)^{(M+1)N}r_{ij} \in \tilde a\Lambda[z_1,z_1^{-1},\ldots,z_d,z_d^{-1}]$$
for each pair $i,j$.
Consider the commutative diagram
$$\matrix 
\noalign{\vskip6pt}
\Lambda\otimes_{\Z_p}\O_{\hat{D}(a)}^m&\buildrel{(F^*\tilde a/\tilde a)^{(M+1)N}(r_{ij})\circ F^*}\over\lolrar&
 \Lambda\otimes_{\Z_p}\O_{\hat{D}(a)}^m  \\ \noalign{\vskip4pt}
\Big\vert{\scriptstyle\tilde a^{(M+1)N}}&&\Big\vert{\scriptstyle\tilde a^{(M+1)N}}\\
\noalign{\vskip-8pt}
\hskip3pt\downarrow\phantom{{\scriptstyle\tilde a^{(M+1)N}}}&&\hskip3pt\downarrow\phantom{\scriptstyle\tilde a^{(M+1)N}}\\
\noalign{\vskip4pt}
\Lambda\otimes_{\Z_p}\O_{\hat{D}(a)}^m&\buildrel{(r_{ij})\circ F^*}\over\lolrar&
 \Lambda\otimes_{\Z_p}\O_{\hat{D}(a)}^m\endmatrix$$  
in which the horizontal arrows are $F^*$-linear maps 
and the vertical maps are $\Lambda\otimes_{\Z_p}\O_{\hat{D}(a)}$-linear isomorphisms
(since they are simply multiplication by $\tilde a^{(M+1)N}$, a unit
on $\hat{D}(a)$).
This shows that we may choose a basis for $\E$ with respect
to which the matrix $(r_{ij})$ of $\Phi$ has entries in
$\tilde a\Lambda[z_1,z_1^{-1},\ldots,z_d,z_d^{-1}],$ and
from now on we assume  such a choice made,
so that
$\L$ is described by the short exact sequence
$$  0 \lrar \L \lrar  
\Lambda\otimes_{\Z_p}\O_{\hat D(a)}^m \buildrel{1-(r_{ij}(z))\circ F^*}\over\lllrar  
\Lambda\otimes_{\Z_p}\O_{\hat D(a)}^m  \llrar   0, $$
in which $(r_{ij}(z))$ is a matrix of of polynomials vanishing 
at the points of $Z(a)$.

Now as in (4.7) we may lift the polynomials $r_{ij}$, this time to elements
$\tilde{r}_{ij}$ of
$\tilde{a}\tilde \Lambda[z_1,z_1^{-1},\ldots,z_d,z_d^{-1}],$ to obtain a lift $\tilde{\E}$
of $\E$ and a lift $\tilde{\L}$ of $\L$. Let $s$ be an integer chosen 
so that each $(z_1\cdots z_d)^s\tilde{r}_{ij}$ belongs to
${\tilde\Lambda}[z_1,\ldots,z_d],$  and choose a second integer $t$ so that $(p-1)t>s$.
Then the diagram 
$$\matrix \noalign{\vskip4pt}
{\tilde\Lambda}\otimes_{\Z_p}\O_{\Gmhat^d}^m&
\buildrel{ (z_1\cdots z_d)^{(p-1)t}(\tilde{r}_{ij})\circ F^*}\over\lolrar&
{\tilde\Lambda}\otimes_{\Z_p}\O_{\Gmhat^d}^m \\ \noalign{\vskip4pt}
\quad\Big\vert{\scriptstyle(z_1\cdots z_d)^t}&&\quad\Big\vert{\scriptstyle(z_1\cdots z_d)^t}\\
\noalign{\vskip-8pt}
\quad\hskip3pt\downarrow\phantom{\scriptstyle(z_1\cdots z_d)^t}&&\quad\hskip3pt\downarrow\phantom{\scriptstyle(z_1\cdots
z_d)^t}\\ {\tilde\Lambda}\otimes_{\Z_p}\O_{\Gmhat^d}^m&\buildrel{(\tilde{r}_{ij})\circ F^*}\over\lolrar&
 {\tilde\Lambda}\otimes_{\Z_p}\O_{\Gmhat^d},\\ \noalign{\vskip4pt}\endmatrix $$
in which the horizontal arrows are $F^*$-linear morphisms and the vertical
arrows are ${\tilde\Lambda}\otimes_{\Z_p}\O_{\Gmhat^d}$-linear isomorphisms (since $(z_1\cdots z_d)^t$
is invertible on $\Gmhat^d$), commutes.
The matrix
$(z_1\cdots z_d)^{(p-1)t}(\tilde{r}_{ij})$ consists of elements of  \vglue-9pt
$$(z_1\cdots z_d){\tilde\Lambda}[z_1,\ldots,z_d]\cap
\tilde{a}{\tilde\Lambda}[z_1,z_1^{-1},\ldots,z_d,z_d^{-1}],\eqno(4.9)$$
\vglue6pt\noindent by construction.
Thus we see that we may choose a basis of $\tilde{\E}$
so that the matrix $\tilde{r}_{ij}$ describing the $F^*$-linear
endomorphism $\tilde{\Phi}$
of $\tilde{\E}$ consists of elements of (4.9), and from now on
we assume   this done, with the matrix $\tilde{r}_{ij}$  chosen with respect
to this basis.

Let $g$ denote the structural morphism $g:\Bbb P^d \rightarrow \Spec \Fp,$
let $h$ denote the open immersion of $\Bbb F_p$-schemes $h:D(a) \rightarrow \Bbb P^d$,
underlying the open immersion of $p$-adic formal schemes
$\hat{D}(a)\rightarrow \hat{\Bbb P}^d$, so that $f=gh$, and
let $\O(-1)_{\hat{\Bbb P}^d}$
denote the ideal sheaf of the formal hyperplane at infinity of $\hat{\Bbb P}^d$
(that is, the hyper-plane described by the equation $Z_0=0$).
Let $u$ be an integer chosen so that $(p-1)u$ is greater than the degree of each
of the polynomials $\tilde{r}_{ij}$. 
By virtue of this choice of $u$, the matrix $(\tilde{r}_{ij})$
induces an $F^*$-linear endomorphism
$(\tilde{r}_{ij})\circ F^*$ of ${\tilde\Lambda}\otimes_{\Z_p}\O(-u)_{\hat{\Bbb P}^d}.$
Furthermore, since the $\tilde{r}_{ij}$ are divisible by each $z_i$ as well
as by $\tilde{a},$ and since $(p-1)u$ is in fact greater than (rather than
just equal to) the degree of any of the $\tilde{r}_{ij}$, we see that the \'etale
sheaf of $(\tilde{r}_{ij})\circ F^*$-invariants of 
${\tilde\Lambda}\otimes_{\Z_p}\O(-u)_{\hat{\Bbb P}^d}^m$
is exactly equal to $h_!\tilde \L$. More precisely, write $\tilde \L_n = \tilde \L/p^n\tilde \L,$ 
and $\tilde \Lambda_n = \tilde \Lambda/p^n\tilde \Lambda.$ Then for each $n \geq 1,$ we have 
have a short exact sequence of \'etale sheaves\vglue-9pt
$$
  0 \lrar h_!\tilde \L_n
\lrar {\tilde \Lambda_n}\otimes_{\Z_p}\O(-u)^m_{\hat{\Bbb P}^d}
\buildrel{1-(\tilde{r}_{ij})\circ F^*}\over\lllrar {\tilde\Lambda_n}\otimes_{\Z_p}\O(-u)^m_{\hat{\Bbb P}^d}\lrar 0. 
$$
\vglue6pt\noindent
This is the exact sequence which we will use to compute  
$f_!\tilde{\L}=g_!h_!\tilde{\L},$ and hence to compute its $L$-function, in
order to compare it with the $L$-function of $\tilde{\L}$. (The construction
of this exact sequence is an application of the technique of 
[4, Lemma~4.5, p.~120]. It is greatly generalized in 
[6], [7].)
\bigbreak

(4.10) If we apply $g_*$ to the short exact sequence of \'etale
sheaves constructed above
  and keep  in mind that for coherent sheaves \'etale push-forward agrees with
the coherent push-forward, and that for $i<d,$
$$R^if_*({\tilde\Lambda_n}\otimes_{\Z_p}\O(-u)_{\hat{\Bbb P}^d})=
H^i(\hat{\Bbb P}^d,{\tilde\Lambda_n}\otimes_{\Z_p}\O(-u)_{\hat{\Bbb P}^d})
=0,$$
we obtain the exact sequence
$$\multline 0 \lrar R^df_!\tilde{\L_n}\lrar H^d({\Bbb P}^d,
{\tilde\Lambda_n}\otimes_{\Z_p}\O(-u)_{\hat{\Bbb P}^d}^m)
 \\
\buildrel{1-H^d((\tilde{r}_{ij})\circ F^*)}\over{\lorar}
H^d({\Bbb P}^d,{\tilde\Lambda_n}\otimes_{\Z_p}\O(-u)_{\hat{\Bbb P}^d}^m) \lrar
R^{d+1}f_!\tilde{\L_n}\lrar 0 \endmultline$$
of \'etale sheaves on $\Spec \Fp.$ 
If we take the stalks of this exact sequence over the geometric point
$\Spec \Fpbar$ of $\Spec \Fp,$ and pass to the inverse limit over $n,$ then we obtain the exact sequence
$$\multline  0 \lrar (R^df_!\tilde{\L})_{\Fpbar}\lrar H^d({\Bbb P}^d,
{\tilde\Lambda}\otimes_{\Z_p}\O(-u)_{\hat{\Bbb P}^d}^m)\otimes_{\Z_p}\Wbar
 \\
\buildrel{1-H^d((\tilde{r}_{ij})\circ F^*)\otimes \sigma} \over{\lorar}
H^d({\Bbb P}^d,{\tilde\Lambda}\otimes_{\Z_p}\O(-u)_{\hat{\Bbb P}^d}^m)\otimes_{\Z_p}\Wbar \\
\lrar
(R^{d+1}f_!\tilde{\L})_{\Fpbar}\lrar 0. \endmultline$$
The morphism
$1-H^d((\tilde{r}_{ij})\circ F^*)\otimes \sigma$ is surjective, and so we
see that\break $R^{d+1}f_!\tilde{\L} = 0.$ 

\specialnumber{4.11} 
\proclaim{Lemma} The ratio
$$L(\Spec \Fp,f_!\tilde{\L})/\det_{{\tilde\Lambda}}(1-H^d((\tilde{r}_{ij})\circ F^*)T,
H^d({\Bbb P}^d,{\tilde\Lambda}\otimes_{\Z_p}\O(-u)_{\hat{\Bbb P}^d}^m))^{(-1)^{d+1}}$$ is an element
of $1+\tilde\goth m T{\tilde\Lambda}\langle T\rangle,$ and is a rational function{\rm .}
\endproclaim

\demo{Proof}
Let us denote the finite rank 
free ${\tilde\Lambda}$-module $H^d({\Bbb P}^d,{\tilde\Lambda}\otimes_{\Z_p}\O(-u)_{\hat{\Bbb P}^d}^m)$
by $M$, for simplicity of notation, and let us denote the ${\tilde\Lambda}$-linear
endomorphism $H^d((\tilde{r}_{ij})\circ F^*)$ of $M$ by $\Psi.$ 
A little linear algebra shows that $M$ has a canonical direct sum decomposition
$M=M_{\rm unit}\oplus M_{\rm nil},$ determined by the property that $M_{\rm unit}$ is
the maximal $\Psi$-invariant
${\tilde\Lambda}$-submodule of $M$ on which $\Psi$ acts surjectively (or equivalently,
bijectively), while $M_{\rm nil}$ is the maximal $\Psi$-invariant
${\tilde\Lambda}$-submodule of $M$ on which $\Psi$ acts topologically nilpotently.
Indeed, it is clear that such decomposition exists for $M$ as a $\Z_p$-module. 
However, as it is canonical, $M_{\rm unit}$ and $M_{\rm nil}$ are $\tilde\Lambda$ stable. 
Thus they are projective, and hence free $\tilde \Lambda$-modules.

We extend the endomorphism $\Psi$ to the $\sigma$-linear endomorphism
$\Psi\otimes \sigma$ of $M\otimes_{\Z_p}\Wbar.$ Since $\Psi$ acts topologically
nilpotently on $M_{\rm nil},$ $1-\Psi\otimes \sigma$ acts bijectively
on $M_{\rm nil}\otimes_{\Z_p} \Wbar.$ Thus we may replace the exact sequence
$$0 \rightarrow (R^df_!\tilde{\L})_{\Fpbar}\rightarrow M\otimes_{\Z_p} \bar W
\buildrel 1-\Psi\otimes \sigma \over \longrightarrow M\otimes_{\Z_p} \bar W
\rightarrow 0$$ of (4.10)  by the exact sequence
$$0\rightarrow (R^df_!\tilde{\L})_{\Fpbar}\rightarrow M_{\rm unit}\otimes_{\Z_p} \bar W
\buildrel 1-\Psi\otimes \sigma \over \longrightarrow M_{\rm unit}\otimes_{\Z_p}\bar W
\rightarrow 0.$$ Since $\Psi$ acts bijectively on $M_{\rm unit},$ we see
that $\Psi$ makes $M_{\rm unit}$ into a unit $({\tilde\Lambda},F)$-crystal, and 
that $R^df_!\tilde{\L}$ is the corresponding lisse $p$-adic \'etale sheaf
on $\Spec \Fp.$ Thus
$$\align L(\Spec \Fp,f_!\tilde{\L}) &
= \det_{{\tilde\Lambda}}(1-\phi T,(R^df_!\tilde{\L})_{\Fpbar})^{(-1)^{d+1}}\\
& = \det_{{\tilde\Lambda}}(1-\Psi T, M_{\rm unit})^{(-1)^{d+1}}\\
& = \det_{{\tilde\Lambda}}(1-\Psi T, M)^{(-1)^{d+1}}/\det_{{\tilde\Lambda}}(1-\Psi T,M_{\rm nil})^{(-1)^{d+1}}.\endalign$$
Thus the lemma will be proved once we show that
$\det_{{\tilde\Lambda}}(1-\Psi T,M_{\rm nil})^{(-1)^{d+1}}$ belongs to 
$1+\tilde\goth m{\tilde\Lambda}\langle T\rangle.$ But
this follows from the fact that $\Psi$ acts topologically nilpotently on $M_{\rm nil}$,
and so nilpotently on $M/\tilde\goth m,$
so that the characteristic polynomial of $\Psi$ on $M_{\rm nil}$ satisfies the
congruence
\medbreak
\hfill ${\displaystyle\det_{{\tilde\Lambda}}(1-\Psi T,M_{\rm nil}) \equiv 1 \pmod {\tilde\goth m}.}$\hfill
\enddemo
 
(4.12) We now turn to describing the $L$-function of $\L$
in terms of the matrix $(\tilde{r}_{ij})$.

\specialnumber{4.13}
\proclaim{Lemma} The $L$\/{\rm -}\/function of $\tilde{\L}$ on $D(a)$ is
determined by the formula
$$L(D(a),\tilde{\L}) = \prod_{x\in |D(a)|}
\det_{{\tilde\Lambda}\otimes_{\Z_p}\Wx}(1-((\tilde{r}_{ij}(\tilde{x}))\circ
\sigma)^{d(x)}T^{d(x)},{\tilde\Lambda}\otimes_{\Z_p}\Wx^m)^{-1}.$$
{\rm (}\/Recall from  {\rm (4.3)} that $\tilde{x}$ denotes the Teichm{\rm \"{\it u}}ller lift of
the closed point $x${\rm .)}
\endproclaim

\demo{Proof} This follows from the fact that $\O_{\hat{D}(a)}^m$
equipped with the $F^*$-linear
endomorphism  $(\tilde{r}_{ij})\circ F^*$ is the unit
$F$-crystal $\tilde{\E}$ corresponding to $\tilde{\L}$.
More precisely, the Artin-Schreier short exact sequence
$$ 0 \lrar \L  \lrar \O_{\hat{D}(a)}^m
\buildrel{1-(\tilde{r}_{ij})\circ F^*}\over{\lllrar} \O_{\hat{D}(a)}^m \lrar 0$$
shows that for any point $x$ of $D(a),$ if $\bar x$ denotes an \'etale point lying 
over $x,$ then 
\medbreak
\hfill${\displaystyle\det(1-\phi^{d(x)}T^{d(x)},\L_{\bar x}) = \det(1-((\tilde{r}_{ij}(\tilde{x}))\circ
\sigma)^{d(x)}T^{d(x)},{\tilde\Lambda}\otimes_{\Z_p}W_x^m).}$\hfill
\enddemo

\specialnumber{4.14}
\proclaim{Corollary} The $L$\/{\rm -}\/function of $\tilde{\L}$ on $D(a)$ agrees
with the product
$$\prod_{x\in |\Gm^d|} \det_{{\tilde\Lambda}\otimes_{\Z_p}\Wx}(1-((\tilde{r}_{ij}(\tilde{x}))\circ
\sigma)^{d(x)}T^{d(x)},{\tilde\Lambda}\otimes_{\Z_p}\Wx^m)^{-1}$$ up to multiplication by an element
of $1+ p T{\tilde\Lambda}\langle T \rangle.$
\endproclaim

\demo{Proof} Lemma~4.13 shows that the ratio of $L(D(a),\tilde{\L})$ and
the product in the statement of the corollary is equal to
$$\prod_{x\in |Z(a)|} \det_{{\tilde\Lambda}\otimes_{\Z_p}\Wx}(1-((\tilde{r}_{ij}
(\tilde{x}))\circ
\sigma)^{d(x)}T^{d(x)},{\tilde\Lambda}\otimes_{\Z_p}\Wx^m)^{-1}.$$
Thus it suffices to show that this
product belongs to $1+ p{\tilde\Lambda}\langle T \rangle.$
 
For any positive integer $\delta$, let $|Z(a)|_{\delta}$ denote the (finite!) subset
of $|Z(a)|$ consisting of those closed points $x$ for which $d(x) = \delta.$
By construction, each $\tilde{r}_{ij}$ is an element of 
$\tilde{a}{\tilde\Lambda}[z_1,z_1^{-1},\ldots, z_d,z_d^{-1}].$
Thus for each $x$ in $|Z(a)|_{\delta}$,
$\tilde{r}_{ij}(\tilde{x})$ lies in $p {\tilde\Lambda}\otimes_{\Z_p}\Wx,$
so that the entries of 
$$((\tilde{r}_{ij}(\tilde{x}))\circ \sigma)^{d(x)}
=((\tilde{r}_{ij}(\tilde{x}))\circ \sigma)^{\delta}$$ lie in $p^{\delta}
\tilde\Lambda,$
and the finite product
$$\prod_{x\in |Z(a)|_d} \det_{\tilde\Lambda\otimes_{\Z_p}\Wx}(1-((r_{ij}(\tilde{x}))\circ
\sigma)^{d(x)}T^{d(x)},\tilde\Lambda\otimes_{\Z_p}\Wx^m)^{-1}$$ 
lies in $1+p^{\delta}T\Lambda[[p^{\delta}T]]\subset 1+p^{\delta}T^{\delta}
{\tilde\Lambda}\langle T^{\delta} \rangle.$ 
From this we conclude that
$$\multline \prod_{x\in |Z(a)|} \det_{{\tilde\Lambda}\otimes_{\Z_p}\Wx}(1-((r_{ij}(\tilde{x}))\circ
\sigma)^{d(x)}T^{d(x)},{\tilde\Lambda}\otimes_{\Z_p}\Wx^m)^{-1}\\
=   \prod_{\delta>0}\  \prod_{x\in |Z(a)|_{\delta}}
\det_{{\tilde\Lambda}\otimes_{\Z_p}\Wx}(1-((r_{ij}(\tilde{x}))\circ
\sigma)^{d(x)}T^{d(x)},{\tilde\Lambda}\otimes_{\Z_p}\Wx^m)^{-1}\endmultline$$ 
lies in $1+pT\tilde\Lambda\langle T \rangle,$
thus proving the corollary.
\enddemo 

(4.15) Theorem~1.5 now follows from Lemma~4.11 and Corollary~4.14,
together with the following result, whose proof is the subject of
Section~5.

\specialnumber{4.16}
\proclaim{Proposition} Let $P(T)$ denote $$\det_{{\tilde\Lambda}}
(1-H^d((\tilde{r}_{ij})\circ F^*)T,
H^d({\Bbb P}^d,{\tilde\Lambda}\otimes_{\Z_p}\O(-u)_{\hat{\Bbb P}^d}^m))\in
1+T{\tilde\Lambda}[T].$$ Then
$$\align&
\prod_{x\in |\Gm^d|} \det_{{\tilde\Lambda}\otimes_{\Z_p}\Wx}(1-((\tilde{r}_{ij}(\tilde{x}))\circ
\sigma)^{d(x)}T^{d(x)},{\tilde\Lambda}\otimes_{\Z_p}\Wx^m)^{-1}
\\  
&\hskip.75in=  \prod_{i=0}^d P(p^iT)^{(-1)^{d+1-i}\left(d\atop
i\right)},\endalign$$ and so in particular,
the ratio $$\left(\prod_{x\in |\Gm^d|} 
\det_{{\tilde\Lambda}\otimes_{\Z_p}\Wx}(1-((\tilde{r}_{ij}(\tilde{x}))\circ
\sigma)^{d(x)}T^{d(x)},{\tilde\Lambda}\otimes_{\Z_p}\Wx^m)^{-1}\right)/P(T)^{(-1)^{d+1}}$$
lies in $1+pT{\tilde\Lambda}\langle T \rangle,$ and is a rational function{\rm .} 
\endproclaim 

\section{Proof of Proposition~4.16}

(5.1) In order to prove Proposition~4.16, we begin by
describing the morphism which is Grothendieck-Serre dual to the
map $1-H^d((\tilde{r}_{ij})\circ F^*).$ 
The free $\tilde \Lambda$-module  
$H^d(\Bbb P^d,{\tilde\Lambda}\otimes_{\Z_p}\O(-u)_{\hat{\Bbb P}^d}^m) 
\iso H^d(\Bbb P^d, \O(-u)_{\hat{\Bbb P}^d}^m)\otimes_{\Z_p}{\tilde\Lambda}$ has 
Grothendieck-Serre dual equal to 
$H^0(\Bbb P^d,\Omega^d(u)_{\hat{\Bbb P}^d}^m)\otimes_{\Z_p}{\tilde\Lambda}$.
Let $\C$ denote the Cartier
operator on $\Omega^d_{\hat{\Bbb P}^d}$ corresponding to the lift of Frobenius $F$,
defined by the formula
$$\multline
\C(z_1^{i_1}\cdots z_d^{i_d}\frac{dz_1\wedge \cdots \wedge dz_d}{z_1\cdots z_d})
\\ = \cases z_1^{\frac{i_1}{p}}\cdots z_d^{\frac{i_d}{p}}
\dfrac{dz_1\wedge\cdots\wedge dz_d}{z_1\cdots z_d}
&\text{if }i_j \equiv 0 \pmod p
\text{ for all $j$}\\
0 &\text{if }i_j\not\equiv 0 \pmod p 
\text{ for some $j$.}\endcases\endmultline$$
By virtue of our choice of $u$,
the composite $\C\circ (\tilde{r}_{ij})$ defines an endomorphism
of ${\tilde\Lambda}\otimes_{\Z_p}\Omega^d(u)$:
$$ {\tilde\Lambda}\otimes_{\Z_p}\Omega^d(u) \buildrel{\C\circ (\tilde{r}_{ij})} \over{\lllrar}
{\tilde\Lambda}\otimes_{\Z_p}\Omega^d(u), $$
and the dual of the endomorphism
$H^d((\tilde{r}_{ij}))\circ F^*$ of $H^d(\Bbb P^d,\O(-u)_{\hat{\Bbb P}^d}^m)$
is the endomorphism
$H^0(\C\circ (\tilde{r}_{ij}))$ of 
$H^0(\Bbb P^d,\Omega^d(u)_{\hat{\Bbb P}^d}^m)\otimes_{\Z_p}{\tilde\Lambda}$.
Since a matrix and its adjoint have the same characteristic polynomial
we see that
$$P(T)=\det_{{\tilde\Lambda}}(1-H^0(\C\circ (\tilde{r}_{ij}))T, H^0(\Bbb
P^d,\Omega^d(u)_{\hat{\Bbb P}^d}^m)\otimes_{\Z_p}{\tilde\Lambda}).$$

\medbreak (5.2) The usual formula for the logarithm of a determinant shows
that
$$\log P(T) = - \sum_{n=1}^\infty \trace_{{\tilde\Lambda}}(H^0(\C\circ(\tilde{r}_{ij}))^n,
H^0(\Bbb P^d,\Omega^d(u)_{\hat{\Bbb P}^d}^m)\otimes_{\Z_p}{\tilde\Lambda}) \dfrac{T^n}{n}. \eqno(5.3)$$
For any sections $a$ of ${\tilde\Lambda}\otimes_{\Z_p}\O_{\hat{\Bbb P^d}}$ and 
$\omega$ of $\Omega^d_{\hat{\Bbb
P^d}}$ we have the identity $$ a \C \omega = \C F^*(a) \omega,$$ 
and so for any positive integer $n$ we have
$$(\C\circ (\tilde{r}_{ij}))^n = \C^n \circ \prod_{i=0}^{n-1}
((F^*)^{i}\tilde{r}_{ij}).$$ This allows us to rewrite 
(5.2) in the form
$$\log P(T) = - \sum_{n=1}^\infty \trace_{{\tilde\Lambda}}(H^0(\C^n \circ \prod_{k=0}^{n-1}
(F^*)^{k}(\tilde{r}_{ij})),H^0(\Bbb
P^d,\Omega^d(u)_{\hat{\Bbb P}^d}^m)\otimes_{\Z_p}{\tilde\Lambda}) \dfrac{T^n}{n}.\eqno(5.4)$$

\medbreak (5.5)  Fix a positive integer $n$,
and suppose that $(f_{ij})$ is any $m\times m$ matrix of elements
of $(z_1\cdots z_d){\tilde\Lambda}[z_1,\ldots,z_d],$ such that each of the $f_{ij}$ is of
degree less than $(p^n-1)u$. (For example, $\prod_{k=0}^{n-1}((F^*)^{k}\tilde{r}_{ij})$
is such
a matrix.) We have a morphism
$$ {\tilde\Lambda}\otimes_{\Z_p}\Omega^d(u)_{\hat{\Bbb P}^d}^m \buildrel{\C^n\circ 
(f_{ij})}\over{\lllrar} {\tilde\Lambda}\otimes_{\Z_p}\Omega^d(u)_{\hat{\Bbb P}^d}^m, $$ which induces an endomorphism
$H^0(\C^n \circ 
(f_{ij}))$ of 
$H^0(\Bbb P^d,{\tilde\Lambda}\otimes_{\Z_p}\Omega^d(u)_{\hat{\Bbb P}^d}^m).$ 
We may also form the sum of matrices
$$\sum_{{x\in\Gm^d(\Bbb F_{p^n})}}(f_{ij}(\tilde{x})),$$
which will be an $m\times m$ matrix of elements of $\Z_p$.
In this situation we have the following lemma, which is a slight modification of
Lemma~2 of [5]:  

\specialnumber{5.6}
\proclaim{Lemma} For any matrix $(f_{ij})$ as above there is the following formula\/{\rm :}\/
$$\align & \trace_{{\tilde\Lambda}}(H^0(\C^n \circ 
(f_{ij})),H^0(\Bbb
P^d,{\tilde\Lambda}\otimes_{\Z_p}\Omega^d(u)_{\hat{\Bbb P}^d}^m)) \\ 
&\hskip.5in = \dfrac{1}{(p^n-1)^d}\trace_{\tilde\Lambda}\left(\sum_{{x\in \Gm^d(\Bbb F_{p^n})}}
(f_{ij}(\tilde{x})), {\tilde\Lambda}^m\right).\endalign$$
{\rm (}\/Recall from  {\rm (4.3)} that $\tilde{x}$ denotes the Teichm{\rm \"{\it u}}ller lift
of the $\Bbb F_{p^n}$-valued point $x$ of $\Gm^d${\rm .)}
\endproclaim 

\demo{Proof}
We first reduce the proof of the formula
to the case $m=1,$ as follows. It is clear that the
trace on the left-hand side is given by the formula
$$\multline \trace_{{\tilde\Lambda}}(H^0(\C^n \circ 
(f_{ij})),H^0(\Bbb
P^d,{\tilde\Lambda}\otimes_{\Z_p}\Omega^d(u)_{\hat{\Bbb P}^d}^m)) \\ =
\sum_{i=1}^m \trace_{{\tilde\Lambda}}(H^0(\C^n \circ f_{ii}),H^0(\Bbb
P^d,{\tilde\Lambda}\otimes_{\Z_p}\Omega^d(u)_{\hat{\Bbb P}^d})),\endmultline$$
while the trace on the right-hand side
is given by the formula
$$\dfrac{1}{(p^n-1)^d} \sum_{i=1}^m \sum_{{x\in \Gm^d(\Bbb F_{p^n})}}
f_{ii}(\tilde{x}).$$ Thus if we knew that for any $f$ in $z_1\cdots
z_d{\tilde\Lambda}[z_1,\ldots,z_d]$ of degree less than $(p^n-1)u$ we had the equation
$$\trace_{\tilde\Lambda}(H^0(\C^n \circ f),H^0(\Bbb
P^d,{\tilde\Lambda}\otimes_{\Z_p}\Omega^d(u)_{\hat{\Bbb P}^d})) =
\dfrac{1}{(p^n-1)^d}\sum_{{x\in \Gm^d(\Bbb F_{p^n})}}
f(\tilde{x}),\eqno(5.7)$$ the lemma would follow. This is what we now prove.

The set of all such $f$ form a ${\tilde\Lambda}$-submodule of 
${\tilde\Lambda}[z_1,\ldots,z_d]$, and both
sides of equation (5.7) are ${\tilde\Lambda}$-linear in $f$.
Thus it suffices to prove the formula in the case that $f$ is a single
monomial $$f = z_1^{\alpha_1}\cdots z_d^{\alpha_d},$$ with $$ 0 < \alpha_1,\ldots,
\alpha_d \text{, \hskip4mm} \alpha_1 + \cdots + \alpha_d < (p^n-1)u.\eqno(5.8)$$

The ${\tilde\Lambda}$-module $H^0(\Bbb
P^d,{\tilde\Lambda}\otimes_{\Z_p}\Omega^d(u)_{\hat{\Bbb P}^d})$ has as a basis the differentials
$$z_1^{\beta_1}\cdots z_d^{\beta_d}\dfrac{dz_1\wedge \cdots \wedge dz_d}
{z_1\cdots z_d},$$ with $$0 <  \beta_1,\ldots,\beta_d \text{, \hskip4mm}
\beta_1 + \cdots + \beta_d < u.\eqno(5.9)$$ We see that
$$\align
& \C^n(z_1^{\alpha_1}\cdots z_d^{\alpha_d} \, z_1^{\beta_1}\cdots
z_d^{\beta_d}\dfrac{dz_1\wedge \cdots \wedge dz_d}
{z_1\cdots z_d})\\
& \qquad  = \C^n(z_1^{\alpha_1+\beta_1} \cdots
z_d^{\alpha_d+\beta_d} \dfrac{dz_1\wedge \cdots \wedge dz_d}
{z_1\cdots z_d})\\
& \qquad  = \cases z_1^{\frac{\alpha_1+\beta_1}{p^n}} \cdots \,
z_d^{\frac{\alpha_d+\beta_d}{p^n}} \dfrac{dz_1\wedge \cdots \wedge dz_d}
{z_1\cdots z_d} & \text{ all } \alpha_i + \beta_i\equiv 0 \pmod{p^n} \\
0 & \text{ some } \alpha_i + \beta_i \not\equiv 0 \pmod{p^n}. \endcases 
\endalign
$$ From this we see that the matrix of $H^0(\C^n\circ f)$ with
respect to the given basis of $H^0(\Bbb P^d,{\tilde\Lambda}\otimes_{\Z_p}\Omega^d(u)_{\hat{\Bbb P}^d})$
has entries which are either 0 or 1, and so the trace of $H^0(\C^n\circ f)$
is equal to the number of 1's on the diagonal, which is given by the number
of $d$-tuples $\beta_1,\cdots, \beta_d$ which satisfy
$$ \alpha_i+\beta_i = p^n\beta_i, \text{ for } 1\leq i \leq d,\eqno(5.10)$$ together
with (5.9). We see immediately that the system of inequalities and equations
(5.9), (5.10) has no solutions unless $$\alpha_1\equiv \cdots \equiv \alpha_d
\equiv 0 \pmod{p^n-1},$$ in which case there is a unique solution (here
we are using the fact that $\alpha_1,\ldots, \alpha_d$ satisfies (5.8)),
so that 
$$\trace_{{\tilde\Lambda}}(\C^n\circ f,
H^0(\Bbb P^d,{\tilde\Lambda}\otimes_{\Z_p}\Omega^d(u)_{\hat{\Bbb P}^d})) =
\cases 1 & \text{ all } \alpha_i \equiv 0 \pmod{p^n-1}\\
0 & \text{ some } \alpha_i \not\equiv 0 \pmod{p^n-1}.\endcases\eqno(5.11)$$
On the other hand, as $x$ ranges over all the $\Bbb F_{p^n}$-valued points
of $\Gm^d$, the set of Teichm\"uller lifts $\tilde{x}$
ranges over all $d$-tuples of $(p^n-1)^{\text{st}}$ roots of unity in $W(\Bbb F_{p^n}).$
Thus $$\align  \frac{1}{(p^n-1)^d}\sum_{{x\in \Gm^d(\Bbb F_{p^n})}}f(\tilde{x})
&= \frac{1}{(p^n-1)^d}
\sum_{{\zeta_1,\cdots,\zeta_d\in \mu_{p^n-1}}}\zeta_1^{\alpha_1} \cdots
\zeta_d^{\alpha_d}\tag 5.12 \\ &=\cases
1 & \text{ all } \alpha_i \equiv 0 \pmod{p^n-1}\\
0 & \text{ some } \alpha_i \not\equiv 0 \pmod{p^n-1}.\endcases \endalign 
 $$
Comparing (5.11) and (5.12), we see that (5.7) is proved, and with it the
lemma.
\enddemo

(5.13) Applying Lemma~5.6 to the matrices
$\prod_{k=0}^{n-1}((F^*)^{k}\tilde{r}_{ij})$ which appear in the expression
on the right side of equation (5.4) we obtain the equation
$$\align \log P(T)& = -\sum_{n=1}^{\infty} \trace_{{\tilde\Lambda}}(\sum_{x\in \Gm^d(\Bbb F_{p^n})}
\prod_{k=0}^{n-1}(\tilde{r}_{ij}(\sigma^k(\tilde{x}))), {\tilde\Lambda}^m)\frac{T^n}{n(p^n-1)^d}\\
& = -\sum_{n=1}^{\infty} \trace_{{\tilde\Lambda}}(\sum_{x\in \Gm^d(\Bbb F_{p^n})}  
((\tilde{r}_{ij}(\tilde{x}))\circ \sigma)^n,{\tilde\Lambda}^m)\frac{T^n}{n(p^n-1)^d}.\endalign$$
This in turn implies that
$$\log \prod_{i=0}^d P(p^iT)^{(-1)^{d+1-i}\left(d\atop
i\right)}= 
\sum_{n=1}^{\infty} \trace_{{\tilde\Lambda}}(\sum_{x\in \Gm^d(\Bbb F_{p^n})}
((\tilde{r}_{ij}(\tilde{x}))\circ \sigma)^n,{\tilde\Lambda}^m)\frac{T^n}{n}.$$
Exponentiating both sides yields Proposition~4.16.

\section{Lifting representations of arithmetic fundamental groups}

(6.1) In this section we discuss some applications of our method 
of lifting locally constant \'etale sheaves of finite free $\Lambda$ 
modules by lifting the associated $(\Lambda,F)$-crystals. For the duration of
this section $\Lambda$ will denote an artinian local $\Z_p$-algebra with
finite residue field and
$X$ will denote a smooth affine $\F_p$-scheme. 
We have seen above that there exists a surjection 
$\tilde \Lambda \rightarrow \Lambda$ 
with $\tilde\Lambda$ a finite flat local $\Z_p$-algebra.
With this notation, we have the following theorem:

\specialnumber{6.2}
\proclaim {Theorem} Let 
$$\rho: \pi_1(X) \rightarrow {\rm GL}_d(\Lambda)$$ 
be a continuous representation of the arithmetic {\rm \'{\it e}}tale fundamental group of $X.$ 
There exists a continuous lifting of $\rho$ 
$$\tilde \rho: \pi_1(X) \rightarrow {\rm GL}_d(\tilde \Lambda).$$ 
\endproclaim

\demo {Proof} The representation $\rho$ corresponds to a locally constant \'etale sheaf 
of finite free $\Lambda$-modules on $X.$ We denote this sheaf by $\L.$
Now we claim that $X$ can be lifted to a formally smooth $p$-adic formal 
scheme $\hat X$ over $\Z_p$ equipped with a lift $F$ of the absolute Frobenius.
Indeed by [10, III] the obstructions to the existence of such a lifting are 
contained in the cohomology of certain coherent sheaves on $X,$ and so vanish,
since $X$ is affine. 

Fix $\hat X$ and $F,$ and denote by $\E$ the unit $(\Lambda, F)$-crystal 
on $\hat X$ corresponding to $\L,$ so that $\E$ is a finite free 
$\Lambda\otimes_{\Z_p} \O_{\hat X}$ sheaf equipped with an isomorphism $F^*\E \iso \E.$
We have already seen in the proof of Proposition~3.2, and also in (4.7), that 
locally on $\hat X,$ $\E$ can be lifted to a finite free 
$\tilde \Lambda\otimes_{\Z_p}\O_{\hat X}$-module $\tilde \E$ equipped with 
an isomorphism $F^*\tilde\E \iso \tilde \E$ lifting $F^*\E \iso \E.$ The 
obstruction to the existence of such a global lift is contained in the 
cohomology of certain coherent $\O_{\hat X}$-modules 
(because $\Lambda\otimes_{\Z_p}\O_{\hat X}$ is $\O_{\hat X}$ coherent),
and so vanishes, since $\O_{\hat X}$ is formal affine. 

This proves the proposition, because by Proposition~3.2, 
$\tilde \E$ corresponds to a lisse sheaf of $\tilde \Lambda$-modules on $\hat X,$ 
which gives the required representation $\tilde \rho.$
\enddemo

\specialnumber{6.3}
\proclaim {Theorem} With the notation of Theorem~{\rm 6.2,} suppose that 
$X \subset \Bbb A^1$ is an open subset of the affine $\Bbb F_p$\/{\rm -}\/line{\rm .} 
Then $\tilde \rho$ may be chosen so that the corresponding 
lisse sheaf of $\tilde \Lambda$\/{\rm -}\/modules $\L$ has a rational $L$\/{\rm -}\/function{\rm .}
\endproclaim

\demo {Proof} As in Section~4, we let $F$ denote the lift of Frobenius on $\hat{\Bbb P}^1$ given (in
affine coordinates) by $z\mapsto z^p$, and let $\hat X \subset \hat{\Bbb P}^1$ denote
the open formal subscheme
corresponding to $X \subset \Bbb P^1.$
We may extend the unit $(\Lambda, F)$-crystal $\E$ attached to 
$\rho$ to a sheaf $\E^+$ of locally free $\Lambda\otimes_{\Z_p}\O_{\hat{\Bbb P}^1}$-modules 
on $\Bbb P^1.$ Indeed to extend across a point of $\Bbb P^1 - X,$ we may work 
in a neighbourhood of this point, and hence assume that $\E$ is 
free over $\Lambda\otimes_{\Z_p} \O_{\hat X}$ over this neighbourhood, when the existence 
of the extension is clear.
Moreover, if $D$ denotes the divisor on $\Bbb P^1$ 
corresponding to $\Bbb P^1 - X,$ and $\hat D$ denotes the divisor on $ \hat{\Bbb P}^1$ 
corresponding to the Teichm\"uller lift of $D,$ then replacing $\E^+$ by $\E^+(-n\hat D),$ 
for some large integer $n,$ we may 
assume that the  isomorphism $F^*\E \iso \E$ extends to a map $F^*\E^+ \rightarrow \E^+$ 
(which cannot be an isomorphism unless $\rho$ is trivial). 
Now set 
$$M = \underline{\rm Hom}_{\Lambda\otimes_{\Z_p}\O_{\hat{\Bbb P}^1}}(F^*\E^+, \E^+).$$
Replacing $\E^+$ by $\E^+(-n\hat D)$ for some large integer $n$ replaces $M$ by\break
$M((p-1)n\hat{D})$, and so we may assume that 
$H^1(\hat{\Bbb P}^1, M) = 0.$ (This will be used later.)

Now we claim that $\E^+$ can be lifted to a locally free 
$\tilde \Lambda \otimes_{\Z_p} \O_{\hat{\Bbb P}^1}$-module $\tilde \E^+.$ 
Indeed, such a lifting exists locally, and the obstruction to a global lifting 
is contained in $H^2$ of certain coherent sheaves on $\Bbb P^1,$ hence vanishes 
(cf. [10, III, 7.1]). 

Let $\goth n$ be the kernel 
of $\tilde \Lambda \rightarrow \Lambda.$ We claim that $F^*\E^+ \iso \E^+$ 
lifts to a morphism $F^*\tilde\E^+/\goth n^2 \rightarrow \tilde\E^+/\goth n^2$ 
which is an isomorphism when restricted to $X$. To see this, we begin by observing
that such 
a lift exists locally. 
More precisely, such a lift exists over any open $U  \subset \Bbb P^1$ which has the
property that
$\tilde \E^+/\goth n^2$ becomes a free $\tilde\Lambda\otimes_{\Z_p}\O_{\hat{\Bbb P}^1}$-module 
when restricted to $U.$ Now let $U_0, U_1, \dots U_r$ be a collection of open subsets of 
$\Bbb P^1$ such that the $U_i$ cover $\Bbb P^1$ and such that
$\tilde \E^+/\goth n^2$ becomes a free $\tilde\Lambda\otimes_{\Z_p}\O_{\hat{\Bbb P}^1}$-module 
when restricted to each $U_i.$ 
Let $\phi_i: F^*\tilde \E^+ \rightarrow \tilde\E^+$ 
be a lift of  $F^*\E^+ \iso \E^+$ over each $U_i$.
The elements $\phi_i - \phi_j$ 
on the intersections $U_i \cap U_j$ give rise to an element of 
$H^1(\Bbb P^1, N),$ where 
$ N = \underline{\rm Hom}_{\Lambda\otimes_{\Z_p}\O_{\hat{\Bbb P}^1}}
(F^*\E^+,\goth n\tilde\E^+/\goth n^2\tilde\E^+).$ 
We claim that $H^1(\Bbb P^1, N) = 0.$ Indeed, choosing a set of generators 
$n_1, \dots, n_s \in \tilde \Lambda$ for the ideal $\goth n,$ we can define a 
surjection of $\Lambda\otimes_{\Z_p}\O_{\hat{\Bbb P}^1}$-modules  
$$ 
 (\E^+)^s\  \raise6pt\hbox{$\underline{\ \scriptstyle (e_1, \dots, e_s) \mapsto e_1n_1 + \dots +
e_sn_s\ }$}\hskip-6pt \rightarrow \ 
\goth n\E^+/\goth n^2\E^+ .
$$ 
Since $F^*\E^+$ is a locally free $\Lambda\otimes_{\Z_p}\hat{\Bbb P}^1$-module, the above surjection 
yields a surjection $M^s \rightarrow N.$ Taking cohomology and recalling that $H^2$ of
a coherent sheaf of $\O_{\hat{\Bbb P}^1}$-modules vanishes, we see
that $H^1(\Bbb P^1, N)$ is a quotient of $H^1(\Bbb P^1, M)^s,$ which vanishes by assumption.
The vanishing of $H^1(\Bbb P^1,N)$ implies that we can modify the $\phi_i$ 
so that they agree on overlaps. 
This gives the 
required lift, since over each of the $U_i$ the $\phi_i$ are automatically 
isomorphisms by Nakayama's lemma, as this is true modulo $\goth n.$ 

Let $M_2 = \underline{\rm Hom}_{\tilde \Lambda\otimes_{\Z_p}\O_{\hat{\Bbb P}^1}}
(F^*\tilde\E^+/\goth n^2,\tilde\E^+/\goth n^2).$ 
We also have $H^1(\hat{\Bbb P}^1, M_2) = 0,$ as we   take the long exact cohomology
sequence attached to the short exact sequence
$0\rightarrow N \rightarrow M_2 \rightarrow M \rightarrow 0,$ and recall  that
$H^1(\Bbb P^1,M) = H^1(\Bbb P^1,N)\break=0.$

Thus we can replace $\Lambda$ by $\tilde \Lambda/\goth n^2$ and $\E^+$ 
by $\tilde\E^+/\goth n^2.$ Repeating the above arguments and taking an inverse
limit we obtain a 
lifting $F^*\tilde \E^+ \rightarrow \tilde \E^+$ of $F^*\E^+ \rightarrow \E^+$ 
which is an isomorphism when restricted to $X$.
Now this restriction $\tilde \E = \tilde \E^+|_{X}$ 
gives rise to a lisse \'etale sheaf $\tilde \L$ of finite free $\tilde \Lambda$-modules
on $X$
which induces the required representation $\tilde \rho.$ To see that the 
$L$-function of $\tilde \rho$ has the required property, first note 
we may replace $\hat X$ by a nonempty open subset, and $\L$ by its restriction to this 
open subset, as this changes the $L$-function only by a rational function. Thus we 
may assume that $X \subset \Bbb G_m,$ and that
$\E$ is free over $\Lambda\otimes_{\Z_p}\O_{\hat X}$ 
(then $\tilde \E$ is free over $\tilde \Lambda\otimes_{\Z_p}\O_{\hat X}$, as one can lift a basis 
for $\E.$) Now we are in exactly the situation of (4.8) (note that Grothendieck's algebraization
theorem guarantees that the matrix entries describing the morphism 
$F^*\tilde{\E} \rightarrow \tilde{\E}$
are rational functions, since by construction this extends to a morphism of coherent
sheaves $F^*\tilde{\E^+} \rightarrow \tilde \E^+$) and the argument of Section~4 shows that 
$L(X, \tilde \L)$ is a rational function over $\tilde \Lambda.$ Indeed in the calculation of 
this $L$-function in Section~4 the only possible nonrational contribution comes from Corollary~4.14,
where (in the notation of the proof of that corollary) one encounters the 
product 
$\prod_{x\in Z(a)} \det_{\tilde{\Lambda}\otimes_{\Z_p}W_x}(1-((\tilde{r}_{ij}(\tilde{x}))
\circ \sigma)^{d(x)}T^{d(x)},\tilde{\Lambda}\otimes_{\Z_p}W_x^m)^{-1}.$ In the general context
of that corollary, this product could be infinite. However, in
the current situation $Z(a)$ is finite (being a proper closed subset of a curve),
and so this is a finite product, and
hence is indeed a rational function!
\enddemo
\references

1   
\name{M.\ Artin, A.\ Grothendieck}, and \name{J.\ L.\ Verdier}, {\it Th{\rm \'{\it e}}orie des Topos et
Cohomologie {\rm \'{\it e}}tale des Sch{\rm \'{\it e}}mas. Tome\/} 3, {\it Lecture Notes in
Math\/}.\ {\bf 305}, Springer-Verlag, New York, 1973.

2
\name{P.\ Berthelot} and \name{W.\ Messing}, Th\'eorie de Dieudonn\'e cristalline.
III. Th\'eor\`emes\break d'\'equivalence et de pleine fid\'elit\'e, in {\it
The Grothendieck Festschrift\/}, {\it Vol\/}.\ I, {\it Progr.\ in
Math\/}.\ {\bf 86}, Birkh\"auser, Boston, MA, 1990, 173--247.

3
\name{R.\ Crew}, $L$-functions of $p$-adic characters and geometric
Iwasawa theory, {\it Invent.\ Math\/}.\ {\bf 88} (1987), 395--403.

4
\name{P.\ Deligne}, {\it S{\rm \'{\it e}}minaire de G{\rm \'{\it e}}ometrie Alg{\rm \'{\it e}}brique du
Bois\/{\rm -}\/Marie\/} SGA $4\frac 12$, {\it Lecture Notes in Math\/}.\ {\bf
569},
Springer-Verlag, New York, 1977.

5 \name{B.\ Dwork}, On the rationality of the zeta function of an
algebraic variety, {\it Amer.\ J.\ Math\/}.\ {\bf 82} (1960), 631--648.

6
 \name{M.\ Emerton} and \name{M.\ Kisin}, A Riemann-Hilbert correspondence
for unit $F$-crystals. I, in preparation.

7 
\bibline, A Riemann-Hilbert correspondence
for unit $F$-crystals. II, in preparation.
 
8 
\name{J.-Y.\ \'Etesse} and \name{B.\ Le Stum}, Fonctions $L$ associ\'ees aux
$F$-isocristaux surconvergents II: Z\'eros et p\^oles unit\'es, {\it
Invent.\ Math\/}.\ 
{\bf 127} (1997), 1--31.

9
\name{T.\ Ekedahl}, On the adic formalism, in {\it The
Grothendieck Festschrift}, {\it Vol.\/} II, {\it Progr.\ in Math\/}.\
{\bf 87}, Birkh\"auser Boston, MA, 1990, 197--218.

10 
\name{A.\ Grothendieck}, {\it Rev\^etments \'Etales et Groupe
Fondamental\/}, {\it Lecture Notes in Math\/}.\ {\bf 224},
Springer-Verlag, New York, 1970.

11
 \name{N.\ Katz}, Travaux de Dwork, {\it S\'eminaire Bourbaki\/},
24 {\it eme ann\'ee\/} (1971/1972), {\it Lecture Notes in Math\/}.\ {\bf
317}
Springer-Verlag, New York, 1973, 167--200.

12 
\bibline, $p$-adic properties of modular schemes and modular
forms, in {\it Modular Functions of One Variable\/} III ({\it Proc.\
Internat.\ Summer School\/},  {\it Univ.\ Antwerp\/}, {\it Antwerp\/}, 1972),
{\it Lecture Notes in Math\/}.\ {\bf 350}, Springer-Verlag, New York,
1973, 69--190.

13 
\name{D.\ Wan}, Meromorphic continuation of $L$-functions of
$p$-adic
representations, {\it Ann.\ of Math\/}.\ {\bf 143} (1996), 469--498.

14
 \bibline, An embedding approach to Dwork's conjecture,
preprint.
\endreferences

\references

SGAfour
\name{M.~Artin}, \name{A.~Grothendieck} and \name{J.~L.~Verdier},
{\it Th\'eorie des topos et cohomologie \'etale des
sch\'emas.~Tome 3}, Lecture Notes in Math., vol.~{\bf 305},
Springer-Verlag, 1973.

BM
\name{P.~Berthelot} and \name{W.~Messing},
Th\'eorie de Dieudonn\'e cristalline.~III.
Th\'eor\`emes d'\'equi\-valence
et de pleine fid\'elit\'e,
{\it The Grothendieck Festschrift, Vol.~I},
Progress in Math., vol.~{\bf 86},
Birkh\"auser, 1990, 173--247.

Cr
\name{R.~Crew},
$L$-functions of $p$-adic characters and
geometric Iwasawa theory,
Inv.~Math.~88 (1987), 395--403.

De
\name{P.~Deligne}, 
{\it S\' eminaire de G\' eom\'etrie Alg\' ebrique
du Bois-Marie SGA $4\frac{1}{2}$},
Lecture Notes in Math., vol.~{\bf 569}, Springer-Verlag,
1977.

Dw
\name{B.~Dwork},
On the rationality of the zeta function of an algebraic variety,
Am.~J.~Math.~82 (1960),
631--648.

EKone
\name{M.~Emerton} and \name{M.~Kisin},
A Riemann-Hilbert correspondence for unit $F$-crystals.~I,
in preparation.

EKtwo
\name{M.~Emerton} and \name{M.~Kisin},
A Riemann-Hilbert correspondence for unit $F$-crystals.~II,
in preparation.

ELS
\name{J.-Y.~Etesse} and \name{B.~Le Stum},
Fonctions $L$ associ\'ees aux $F$-isocristaux surconvergents
II: Z\'eros et p\^oles unit\'es,
Inv.~Math.~127 (1997), 
1--31.

Ek
\name{T.~Ekedahl}, 
On the adic formalism,
{\it The Grothendieck Festschrift, Vol.~II},
Progress in Math., vol.~{\bf 87},
Birkh\"auser, 1990, 197--218.

SGAone
\name{A.~Grothendieck},
{\it Rev\^etments \'Etales et Groupe Fondamental},
Lecutre Notes in Math., vol.~{\bf 224},
Springer-Verlag, 1970.

Ka
\name{N.~Katz}, 
Travaux de Dwork,
{\it Seminaire Bourbaki, 24\`eme ann\'ee (1971/1972)},
Lecture Notes in Math., vol.~{\bf 317},
Springer-Verlag, 1973, 69--190.

Katwo
\name{N.~Katz},
$p$-adic properties of modular schemes, 
and modular forms,
{\it Modular functions of one variable III
(Proc.~Internat.~Summer School, Univ.~Antwerp, Antwerp 1972)},
Lecture Notes in Math., vol.~{\bf 350},
Springer-Verlag, 
1973, 69--190.

Wa
\name{D.~Wan},
Meromorphic continuation of $L$-functions of $p$-adic 
representations,
Ann.~Math. 143 (1996),
469--498.

Watwo
\name{D.~Wan},
An Embedding Approach to Dwork's Conjecture,
preprint.

\endreferences

\bye